\documentclass{article}

\usepackage{hyperref}

\usepackage{hyperref}

\usepackage{amssymb}

\newcommand{\newc}{\newcommand}

\newc{\eqnoset}{\setcounter{equation}{0}}

\newcommand{\mref}[1]{(\ref{#1})}
\newcommand{\reflemm}[1]{Lemma~\ref{#1}}
\newcommand{\refrem}[1]{Remark~\ref{#1}}
\newcommand{\reftheo}[1]{Theorem~\ref{#1}}

\newcommand{\refdef}[1]{Definition~\ref{#1}}
\newcommand{\refcoro}[1]{Corollary~\ref{#1}}

\newcommand{\refsec}[1]{Section~\ref{#1}}

\newcommand{\beq}{\begin{equation}}
\newcommand{\eeq}{\end{equation}}
\newcommand{\beqno}[1]{\begin{equation}\label{#1}}

\newcommand{\barr}{\begin{array}}
\newcommand{\earr}{\end{array}}

\newc{\bearr}{\begin{eqnarray*}}
\newc{\eearr}{\end{eqnarray*}}

\newc{\bearrno}[1]{\begin{eqnarray}\label{#1}}
\newc{\eearrno}{\end{eqnarray}}

\newc{\non}{\nonumber}
\newc{\nol}{\nonumber\nl}

\newcommand{\bdes}{\begin{description}}
\newcommand{\edes}{\end{description}}
\newc{\benu}{\begin{enumerate}}
\newc{\eenu}{\end{enumerate}}
\newc{\btab}{\begin{tabular}}
\newc{\etab}{\end{tabular}}

\newtheorem{theorem}{Theorem}[section]
\newtheorem{defi}[theorem]{Definition}
\newtheorem{lemma}[theorem]{Lemma}
\newtheorem{rem}[theorem]{Remark}
\newtheorem{exam}[theorem]{Example}
\newtheorem{propo}[theorem]{Proposition}
\newtheorem{corol}[theorem]{Corollary}
\newtheorem{conj}[theorem]{Conjecture}

\newcommand{\btheo}[1]{\begin{theorem}\label{#1}}
\newc{\brem}[1]{\begin{rem}\label{#1}\em}
\newc{\bexam}[1]{\begin{exam}\label{#1}\em}
\newc{\bdefi}[1]{\begin{defi}\label{#1}}
\newcommand{\blemm}[1]{\begin{lemma}\label{#1}}
\newcommand{\bprop}[1]{\begin{propo}\label{#1}}
\newcommand{\bcoro}[1]{\begin{corol}\label{#1}}
\newcommand{\btheoc}[1]{\begin{conj}\label{#1}}
\newcommand{\etheo}{\end{theorem}}
\newc{\etheoc}{\end{conj}}
\newcommand{\elemm}{\end{lemma}}
\newcommand{\eprop}{\end{propo}}
\newcommand{\ecoro}{\end{corol}}
\newc{\erem}{\end{rem}}
\newc{\eexam}{\end{exam}}
\newc{\edefi}{\end{defi}}

\newc{\rmk}[1]{{\bf REMARK #1: }}
\newc{\DN}[1]{{\bf DEFINITION #1: }}

\newcommand{\bproof}{{\bf Proof:~~}}
\newc{\eproof}{{\vrule height8pt width5pt depth0pt}\vspace{3mm}}

\newc{\bfrac}[2]{\dspl{\frac{#1}{#2}}}

\newc{\nid}{\noindent}

\newcommand{\dspl}{\displaystyle}
\newc{\grad}{\nabla}
\newc{\Div}{\mbox{div}}
\newc{\pdt}[1]{\dspl{\frac{\partial{#1}}{\partial t}}}
\newc{\pdn}[1]{\dspl{\frac{\partial{#1}}{\partial \nu}}}
\newc{\pdNi}[1]{\dspl{\frac{\partial{#1}}{\partial \mathcal{N}_i}}}
\newc{\pD}[2]{\dspl{\frac{\partial{#1}}{\partial #2}}}
\newc{\dt}{\dspl{\frac{d}{dt}}}
\newc{\bdry}[1]{\mbox{$\partial #1$}}
\newc{\sgn}{\mbox{sign}}

\newc{\Hess}[1]{\frac{\partial^2 #1}{\pdh z_i \pdh z_j}}
\newc{\hess}[1]{\partial^2 #1/\pdh z_i \pdh z_j}

\newc{\ag}{\alpha}
\newc{\bg}{\beta}
\newc{\cg}{\gamma}\newc{\Cg}{\Gamma}
\newc{\dg}{\delta}\newc{\Dg}{\Delta}
\newc{\eg}{\varepsilon}
\newc{\zg}{\zeta}
\newc{\thg}{\theta}
\newc{\llg}{\lambda}\newc{\LLg}{\Lambda}
\newc{\kg}{\kappa}
\newc{\rg}{\rho}
\newc{\sg}{\sigma}\newc{\Sg}{\Sigma}
\newc{\tg}{\tau}
\newc{\fg}{\phi}\newc{\Fg}{\Phi}
\newc{\vfg}{\varphi}
\newc{\og}{\omega}\newc{\Og}{\Omega}
\newc{\pdh}{\partial}

\newc{\ccG}{{\cal G}}

\newc{\ii}[1]{\int_{#1}}
\newc{\iidx}[2]{{\dspl\int_{#1}~#2~dx}}
\newc{\bii}[1]{{\dspl \ii{#1} }}
\newc{\biii}[2]{{\dspl \iii{#1}{#2} }}
\newc{\su}[2]{\sum_{#1}^{#2}}
\newc{\bsu}[2]{{\dspl \su{#1}{#2} }}

\newc{\biiom}[1]{{\dspl\int_{\bdrom}~ #1 ~d\sg}}
\newc{\io}[1]{{\dspl\int_{\Og}~ #1 ~dx}}
\newc{\bio}[1]{{\dspl\int_{\bdrom}~ #1 ~d\sg}}
\newc{\bsir}{\bsu{i=1}{r}}
\newc{\bsim}{\bsu{i=1}{m}}

\newc{\iibr}[2]{\iidx{\bprw{#1}}{#2}}
\newc{\Intbr}[1]{\iibr{R}{#1}}
\newc{\intbr}[1]{\iibr{\rg}{#1}}
\newc{\intt}[3]{\int_{#1}^{#2}\int_\Og~#3~dxdt}
\newc{\itQ}[2]{\dspl{\int\hspace{-2.5mm}\int_{#1}~#2~dz}}
\newc{\mitQ}[2]{\dspl{\rule[1mm]{4mm}{.3mm}\hspace{-5.3mm}\int\hspace{-2.5mm}\int_{#1}~#2~dz}}
\newc{\mitQQ}[3]{\dspl{\rule[1mm]{4mm}{.3mm}\hspace{-5.3mm}\int\hspace{-2.5mm}\int_{#1}~#2~#3}}

\newc{\mitx}[2]{\dspl{\rule[1mm]{3mm}{.3mm}\hspace{-4mm}\int_{#1}~#2~dx}}
\newc{\mitmu}[2]{\dspl{\rule[1mm]{3mm}{.3mm}\hspace{-4mm}\int_{#1}~#2~d\mu}}
\newc{\iidmu}[2]{\iidx{#1}{#2}}

\newc{\iidm}[3]{{\dspl\int_{#1}~#2~d #3}}

\newc{\itQmu}[2]{\dspl{\int\hspace{-2.5mm}\int_{#1}~#2~d\mu}}
\newc{\mitQmu}[2]{\dspl{\rule[1mm]{4mm}{.3mm}\hspace{-5.3mm}\int\hspace{-2.5mm}\int_{#1}~#2~d\mu}}

\newc{\mitQq}[2]{\dspl{\rule[1mm]{4mm}{.3mm}\hspace{-5.3mm}\int\hspace{-2.5mm}\int_{#1}~#2~d\bar{z}}}
\newc{\itQq}[2]{\dspl{\int\hspace{-2.5mm}\int_{#1}~#2~d\bar{z}}}

\newc{\pder}[2]{\dspl{\frac{\partial #1}{\partial #2}}}

\newc{\bdrom}{\bdry{\Og}}

\newc{\bilhom}{\mbox{Bil}(\mbox{Hom}(\RR^{nm},\RR^{nm}))}
\newc{\VV}[1]{{V(Q_{#1})}}

\newc{\ccA}{{\mathcal A}}
\newc{\ccB}{{\mathcal B}}
\newc{\ccC}{{\mathcal C}}
\newc{\ccD}{{\mathcal D}}
\newc{\ccE}{{\mathcal E}}
\newc{\ccH}{\mathcal{H}}
\newc{\ccF}{\mathcal{F}}
\newc{\ccI}{{\mathcal I}}
\newc{\ccJ}{{\mathcal J}}
\newc{\ccK}{{\mathcal K}}
\newc{\ccP}{{\mathcal P}}
\newc{\ccQ}{{\mathcal Q}}
\newc{\ccR}{{\mathcal R}}
\newc{\ccS}{{\mathcal S}}
\newc{\ccT}{{\mathcal T}}
\newc{\ccX}{{\mathcal X}}
\newc{\ccY}{{\mathcal Y}}
\newc{\ccZ}{{\mathcal Z}}

\newc{\bb}[1]{{\mathbf #1}}

\newc{\myprod}[1]{\langle #1 \rangle}
\newc{\mypar}[1]{\left( #1 \right)}
\newc{\BLLg}{\mathbf{\LLg}}

\newc{\mA}{\mathbf{A}}
\newc{\mB}{\mathbf{B}}
\newc{\mC}{\mathbf{C}}
\newc{\mD}{\mathbf{D}}
\newc{\mE}{\mathbf{E}}
\newc{\mF}{\mathbf{F}}
\newc{\mJ}{\mathbf{J}}
\newc{\mG}{\mathbf{G}}
\newc{\mP}{\mathbf{P}}
\newc{\mR}{\mathbf{R}}
\newc{\mQ}{\mathbf{Q}}
\newc{\mX}{\mathbf{X}}
\newc{\muu}{\mathbf{u}}
\newc{\mvv}{\mathbf{v}}

\newc{\mllg}{\mathbb{\lambda}}
\newc{\mLLg}{\mathbf{\LLg}}

\newc{\lspn}[2]{\mbox{$\| #1\|_{\Lsp{#2}}$}}
\newc{\Lpn}[2]{\mbox{$\| #1\|_{#2}$}}
\newc{\Hn}[1]{\mbox{$\| #1\|_{H^1(\Og)}$}}

\newc{\mynorm}[2]{\| #1\|_{#2}}

\newcommand{\RR}{{\rm I\kern -1.6pt{\rm R}}}

\newc{\itQQ}[2]{\dspl{\int_{#1}#2\,dz}}
\newc{\mmitQQ}[2]{\dspl{\rule[1mm]{4mm}{.3mm}\hspace{-4.3mm}\int_{#1}~#2~dz}}
\newc{\MmitQQ}[2]{\dspl{\rule[1mm]{4mm}{.3mm}\hspace{-4.3mm}\int_{#1}~#2~d\mu}}

\newc{\MUmitQQ}[3]{\dspl{\rule[1mm]{4mm}{.3mm}\hspace{-4.3mm}\int_{#1}~#2~d#3}}
\newc{\MUitQQ}[3]{\dspl{\int_{#1}~#2~d#3}}

\newc{\mccP}{\mathbb{P}}
\newc{\mccK}{\mathbb{K}}

\newc{\DKTmU}{\mccK(U)}
\newc{\DKTmUold}{(K_U(U)^{-1})^T}

\newc{\myPi}{\mathbf{W}}
\newc{\myIbar}{\bar{\ccI}_1}
\newc{\myIhat}{\hat{\ccI}_1}
\newc{\myIbreve}{\breve{\ccI}_0}

\newc{\mmk}{\mathbf{k}}

\newcommand{\ma}{\mathbf{a}}
\newcommand{\mg}{\mathbf{g}}

\newc{\mfu}{\mathbf{f_u}}
\newc{\mh}{\mathbf{h}}
\newc{\mb}{\mathbf{b}}
\newc{\mf}{\mathbf{f}}

\newcommand{\barrl}[2]{\barr{ll}\lefteqn{#1}\hspace{#2}&\\}

\newc{\twomatrix}[1]{\left[\barr{cc}#1\earr\right]}
\newc{\threematrix}[1]{\left[\barr{ccc}#1\earr\right]}

\newc{\mN}{\mathbf{N}}
\newc{\mI}{\mathbf{I}}
\newc{\mH}{\mathbf{H}}

\newc{\mk}{\mathbf{k}}
\newc{\mr}{\mathbf{r}}

\newc{\DIAGM}[2]{\left[\barr{ccc}#1&0\ldots&0\\
	\vdots&\ddots&\vdots\\0&\ldots0&#2\earr \right]}
\newc{\DiagM}[2]{\mbox{diag}\left[#1
	\cdots #2 \right]}
\newc{\vVEC}[2]{\left[\barr{c}#1\\
	\vdots\\#2\earr \right]}
\newc{\hVEC}[2]{\left[#1
	\cdots #2 \right]}

\newc{\mq}{\mathbf{q}}

\newc{\msys}[1]{\left\{\barr{l}#1\earr
	\right.}
\newc{\msysa}[1]{\left\{\barr{ll}#1\earr
	\right.}

\newc{\bbM}{\mathbb{M}}
\newc{\mat}[1]{\left[\barr{cc}#1\earr\right]}

\newc{\me}{\mathbf{e}}

\newc{\vecc}[2]{\left[\barr{cc}#1\\#2\earr\right]}
\newc{\mL}{\mathbb{L}}

\newc{\cO}{{\cal O}}

\newc{\cM}{{\cal M}}

\newc{\myega }{\eg_0(R)}
\newc{\myeg}{\eg_1(\eg_*)}
\newc{\myegp}{\hat{\eg}_1(\eg_*)}

\newc{\diagA}{\mathbb{A}_d}
\newc{\mBB}{\mathbb{B}}
\newc{\MLT}[1]{{\cal M}_{lt}(\Og,#1)}
\newc{\ALT}[1]{{\cal A}_{l}(\Og,#1)}
\newc{\mM}{\mathbb{M}}

\newc{\diag}[1]{\mbox{diag}(#1)}
\newc{\off}[1]{\mbox{offdiag}(#1)}
\newc{\mT}{\mathbb{T}}

\usepackage{amssymb}

\begin{document}

\vspace*{-.8in}
\begin{center} {\LARGE\em Strong Persistence of a Class of Strongly Coupled Parabolic Systems of $m$ Equations.}

 \end{center}

\vspace{.1in}

\begin{center}

{\sc Dung Le}{\footnote {Department of Mathematics, University of
Texas at San
Antonio, One UTSA Circle, San Antonio, TX 78249. {\tt Email: Dung.Le@utsa.edu}\\
{\em
Mathematics Subject Classifications:} 35J70, 35B65, 42B37.
\hfil\break\indent {\em Key words:} Cross diffusion systems,  H\"older
regularity, global existence.}}

\end{center}

\begin{abstract}
	We establish one of the most important assumptions of the strong persistence theory for dynamical systems associated to cross diffusion systems of $m$ equations ($m\ge2$): the stable sets of semi-trivial steady cannot intersect the interior of the positive cone of $C(\Og,\RR^m)$. Many examples will be provide to show the effects of the cross diffusion. \end{abstract}

\section{Introduction}\label{intoMAX}\eqnoset

In this paper, we consider the following parabolic  partial differential systems on a smooth bounded domain $\Og\subset \RR^N$ of $m$ equations ($N,m\ge 2$) for $W=[u,v]^T$ with $u\in\RR^{m_1}, v\in\RR^{m_2}$ and $m=m_1+m_2$ and $m_i\ge1$ 

\beqno{syswhole}\left\{\barr{ll} W_t=\Div({\cal A}(W)DW)+{\cal B}(W)DW+{\cal F}(W,DW)& (x,t)\in Q=\Og\times(0,T),\\u(x,0)=U_0(x)& x\in\Og\\\mbox{$u=0$ on $\partial \Og\times(0,T)$}. &\earr\right.\eeq
Here, in block form
$${\cal A}(W)=\left[\barr{cc}\ma^{11}(W)&\ma^{21}(W)\\\ma^{21}(W)&\ma^{22}(W)\earr\right],\; {\cal B}(W)=\left[\barr{cc}\mb^{11}(W)&\mb^{21}(W)\\\mb^{21}(W)&\mb^{22}(W)\earr\right],\;{\cal F}(W)=\left[\barr{cc}\mg^{11}(W)&0\\0&\mg^{22}(W)\earr\right]. $$
That is,  $\ma^{ii}, \mb^{ii}, \mg^{ii}$ are  square matrices  of size $m_i\times m_i$, for  $i=1,2$,  with sufficiently smooth function entries on $\Og$. As usual, we assume the normal ellipticity (see \cite{Am2} and note that we do not assume that ${\cal A}$ is symmetric): for some positive constants $\llg_0,\LLg_0$
$$\llg_0|\zeta|^2\le \myprod{{\cal A}\zeta,\zeta}\le \LLg_0|\zeta|^2 \quad\forall\zeta\in \RR^m.$$

Coexistence and persistence results for \mref{syswhole} are investigated in \cite[Chapter 6 and 7]{dlebook1}, \cite{dlecoexper} when $m_1=m_2=1$. The general setting can be extended to the case $m_i\ge 1$ as follows: assume that \mref{syswhole} defines a global dynamical system in the positive cone $\RR^m$ and has a semi-trivial steady state solution $W^*=[u^*,0]^T$ with $u^*>0$. This is the case for example
\beqno{Wstar} \left\{\barr{l}-\Div(\ma^{11}(u^*,0)Du^*)=\mb^{11}(u^*,0)Du^*+\mg^{11}(u^*,0)u^*,
\\\ma^{21}(u^*,0)=\mb^{21}(u^*,0).\earr\right.
\eeq
These assumptions and notations will be assumed throughout this paper.

One of the main goals of this work is to show that, under some suitable conditions, any global positive solution of \mref{syswhole} cannot approach this semi-trivial steady state $W^*$ as $t\to\infty$. We will show that the introduction of an appropriate cross diffusion (${\cal A}$ is a full matrix) will make this phenomenon to occur. This is also an important step in establishing the uniform persistence of the dynamical system associated to \mref{syswhole}.

In \refsec{persistence} we briefly recall the abstract theory in \cite{HHZ} to establish the uniform persistence property of \mref{syswhole} by the verification of \reftheo{hhzthm} (\cite[Theorem 4.3]{HHZ}). Proving that \mref{syswhole} defines a global dynamical system and has the positive cone $\RR^m$ as an invariant set is a very hard problem and we will refer the readers to \cite{dlebook1} for recent results of this direction. 
We refer the readers to \cite[Chapter 3 and 4]{dlebook1} (see also \cite{GiaS} for additional details) on how to obtain uniform bound for the derivatives $Du, Dv$ once we obtained that the solution $(u,v)$ is (uniform) H\"older continuous  (see also \cite{dleJMAA}). 
Also, we assume that $C(\Og,\RR^m_+)$ is invariant under the flows of this dynamical system. The chain condition or the dynamics of the flows on the surfaces of $C(\Og,\RR^m_+)$ is also hard for systems of more than 2 equations on higher dimension domains.

Among others,
the m.4) condition (see \refsec{persistence}) that the stable set of the semi trivial steady state $W^*$ does not intersect the interior of $C(\Og,\RR^m_+)$ is one of the most technical and  of interest by itself in other applications. 

In this paper, we will study the effect of cross diffusion on  m.4). Interestingly, we will show that m.4) may not occur for classical reaction-diffusion systems but with an appropriate introduction of cross diffusion we can force m.4) to happen. To the best of our knowledge this result is new and firstly reported here.

The paper is organized as follows.
To begin, we present a technical lemma which will be one of the main vehicles in our argument later if m.4) fails to obtain a contradiction. We then use this lemma to our \mref{syswhole} Our crucial assumption is that the semi-trivial steady state is unstable in its complementary direction ($u$). This is equivalent to the fact that (\cite{dlebook1})  the principal eigenvalue $\tau$ of certain positive operators of the form $L^{-1}M$, where $L$ is an operator associated to a second order elliptic system and $M$ is a matrix, to be strictly greater than 1 or the associated eigenvalue problems (see \mref{per1}) has a positive solutions. This eigenvalue problem is determined by the subsystem for  $v$ in \mref{syswhole} at $W^*$ (i.e. $\ma^{22},\ma^{21},\mb^{22},\mg^{22}$).

In \refsec{tausec} we will take a close look at the eigenvalue problem \mref{per1} and present some necessary conditions  where the crucial assumption $\tau>1$ can be realized in an abstract setting.

If $\ma^{22},\ma^{21}_v, \mb^{22}$ are diagonal then the situation is somewhat simpler but still nontrivial as we can apply the maximum prinples for diagonal system with {\em cooperative} reactions (e.g., see \cite{lopez}). 

Throughout this paper, we adapt the notation $A<B$ (respectively $A\le B$) for two matrices $A=[a_{i,j}], B=[b_{i,j}]$ of the same size in the componentwise sense, i.e. $a_{i,j}<b_{i,j}$ (respectively $a_{i,j}\le b_{i,j}$) for all $i,j$. A matrix $A=[a_{i,j}]$ is said to be cooperative (competitive) if $a_{i,j}>0$ ($a_{i,j}<0$) if $i\ne j$. $A$ is partially competitive if $a_{i,j}<0$ for some $i\ne j$. Of course, $A$ is strictly positive if $a_{i,j}>0$ for all $i,j$.

In \refsec{diagsec}, we will assume that $\ma^{22}(u^*,0), \ma^{21}_v(u^*,0)$ are diagonal (so is $\Div(\ma^{21}_v(u^*,0)Du^*)$; but $\ma^{22}(u,v), \ma^{21}_v(u,v)$ can be non-diagonal when $v\ne0$), so that $L^{-1} $ is diagonal and $(L+kId)^{-1}$ is a positive operator if $k>0$ is sufficiently large. Also by this way, we need only that the off-diagonal of $\mb^{21}_v(W^*)Du^*+\mg^{22}(W^*)$ are {\em nonnegative}. We will see that if $\mg^{22}=[\mg_{ij}]$ with $\mg_{ij}$'s are positively large then this is the case. However,  if $\mg_{ij}$'s are large then it is harder to establish the existence of a global attractor. We will provide several examples and counterexamples \refsec{counterex} and see that even if there is cross diffusion and $\ma^{21}_v(u^*,0)$ is diagonal then $\mg^{22}$ cannot be completely competitive.
These will used later in comparison to the nondiagonal case ($\ma^{22}(u^*,0)$ or $\ma^{21}_v(u^*,0)$ is nondiagonal).

We then continue with the argument in \refsec{nondiagsec} to prove that $\tau>1$ for cross diffusion systems in \refsec{nondiag}. One of the most crucial assumption there is that $L^{-1}M$ is a strongly positive operator and we will check this assumption first. This problem for cross diffusion (non-diagonal case) has been recently investigated in \cite[Theorem 5.14]{MAXP} and we will make use of this in \reftheo{GenMPMatTKRnew}. As we have seen that for diagonal systems, the reaction parts could not be completely competitive to obtain that $\tau>1$. However, if we introduce cross diffusions then we can allow the reaction to be completely competitive but still obtain this goal. Examples will be presented to show the contrast effect to the diagonal case in the previous section. 
We will show that we can have persistence if $\mA$ is full but appropriately designed but no persistence when $\mA$ is replaced by  full/diagonal matrices $A$ while using  the same reaction. Therefore, the structures of $L$ (or ${\cal L}$) in \reftheo{GenMPMatTKRnew} seems to be necessary for persistence.

\section{Persistence}\label{persistence}

To begin, we need the concept of attractors for dissipative dynamical systems. We will discuss it briefly below and refer the readers to another comprehensive literature (e.g. \cite{ET,Temam,ST} and the reference therein) 

Let $(X, d)$ be a metric space and $\{\Fg(t)\}_{t\ge0}$ be a semiflow on $X$. A
subset $A \subset X$ is said to be an attractor for $\{\Fg(t)\}_{t\ge0}$ if $A$ is nonempty, compact, invariant, and
there exists some open neighborhood $U$ of $A$ in $X$ such that $\lim_{t\to\infty} d(\Fg(t)(u), A) = 0$ for
all $\in U$. Here, $d(x, A)$ is the usual Hausdorff distance from $x$ to the set $A$. If $A$ is an
attractor which attracts every point in $X$, $A$ is called global attractor. 

For a nonempty
invariant set $M$, the set $$W^s
(M) := \{x \in X : \lim_{t\to\infty} d(\Fg(t)(x), M) = 0\}$$ is called the
stable set of $M$. A nonempty invariant subset $M$ of $X$ is said to be isolated if it is the
maximal invariant set in some neighborhood of itself.

Let $A$ and $B$ be two isolated invariant sets. $A$ is said to be chained to $B$, denoted by
$A \to B$, if there exists a globally defined trajectory $\fg(t)(x)$, $t\in(-\infty,\infty)$, through some
$x\not\in A\cup B$ whose range has compact closure such that the omega limit set $\og(x) \subset B$ and
the alpha limit set $\ag(x) \subset A$. A finite sequence $\{M_1, M_2, ..., M_k\}$ of isolated invariant sets
is called a \index{chain}chain if $M_1 \to M2 \to \ldots\to M_k$. The chain is called a\index{cycle} cycle if $M_k = M_1$. These concepts were introduced in  \cite{BW}.

Let $X_0 \subset X$ be an open set and $\partial X_0 = X \setminus X_0$. Assume that $X_0$ is positively
invariant. Let $d(x, \partial X_0)$ be the distance from $x$ to $\partial X_0$. $\{\Fg(t)\}_{t\ge0}$ is said to be uniformly
persistent with respect to $(X_0, \partial X_0)$ if there exists $\eta > 0$ such that
$$\liminf_{t\to\infty}d(\Fg(t)(x),\partial X_0)\ge \eta
\quad \mbox{for all } x \in X_0.$$

As in \cite{dleper2}, in order to discuss the uniform persistence property of our parabolic systems we rely on the following result \btheo{hhzthm}(\cite[Theorem 4.3]{HHZ}) Assume that
\bdes \item[C.1)] $\{\Fg(t)\}_{t\ge0}$ has a global attractor $A$;

\item[C.2)] There exists a finite sequence $M = \{M_1, ..., M_k\}$ of pairwise disjoint, compact
and isolated invariant sets in $\partial X_0$ with the following properties:
\bdes
\item[m.1)] $\cup_{x\in\partial X_0} \og(x)\subset \cup_{i=1}^k M_i$,
\item[m.2)] no set of M forms a cycle in $\partial X_0$,
\item[m.3)] $M_i$ is isolated in $X$,
\item[m.4)] $W^s(M_i)\cap X_0=\emptyset$ for each $i=1,\ldots, k$.
\edes

\edes

Then $\{\Fg(t)\}_{t\ge0}$ is uniformly persistent with respect to $(X_0, \partial X_0)$.

\etheo

In this section, we will set
$$X=C_+(\Og)\times C_+(\Og) \mbox{ and } X_0=\{(u,v)\in X\,:\, u(x)>0 \mbox{ and } v(x)>0 \mbox{ for all $x$}\},$$ with the metric $d(U,V)=\|U-V\|_{C(\Og)}$.

One should note that the results here are stronger than those in \cite{dleper1,dleper2}, where we set $X=C^1_+(\Og)\times C^1_+(\Og)$ with the metric $d(U,V)=\|U-V\|_{C^1(\Og)}$.

With the setting of this paper, via the uniform estimates in \cite{dleJMAA,dlebook1},  we see that C.1) is verified. 

We follows \cite{dlebook1} to verify the assumptions m.1)-m.4) of C.2). The boundary of $X_0$ ($u=0$ or $v=0$) is invariant so that m.1) holds. The assumptions m.2) and m.3) involve with the study of dynamics of solutions to the scalar parabolic equation $$w_t=\Div(a(w)Dw)+b(w)Dw+g(w)w$$ on the axis $u=0$ or $v=0$ of $X$. For decades this problem has been extensively studied in literature and we refer the readers to those partially referred in \cite{dleper1, dleper2}. In particular, a very often scenario occurs in applications that we could take $M=\{M_0,M_1,M_2\}$ with $M_0=\{(0,0)\}$ and $M_1,M_2$ are single point sets  $\{(u^*,0)\}$ and $\{(0,v^*)\}$ for some nonnegative $u^*,v^*$ which are solutions to the corresponding elliptic equations. Thus, the sets of semi-trivial solutions $Z_i=M_i$, $i=1,2$, are determined. The assumptions m.2), m.3) are then often verified. Throughout this paper we will assume these settings.

The main difficulty is the justification of m.4) and this is the subject of the next sections. We will prove m.4) by contradiction, under quite general assumptions involving certain eigenvalue problems. Thus, assume that there some positive initial data $(u_0,v_0)\in X_0$ and the corresponding solution $(u,v)$ converges to $M_1$ or $M_2$ in $X$. To fix the idea, let us assume $(u,v)$ converges to $M_1=\{(u^*,0)\}$ in $X$ and we will establish a contradiction if $(u^*,0)$ is $v$-unstable as in \cite{dlejfa} and \cite[Chapter 6]{dlebook1}.

\section{The dynamics of evolutionary solutions near a semi-trivial steady state}\label{dynevol}\eqnoset

We are going to verify the condition m.4) in this section.
Assuming that there is a semi-trivial solution of the considered evolution system converge to $W^*$, we will try to establish a contradiction. The following calculation is verbatim to that of \cite{dlebook1} for systems of two equations to the case of $m=m_1\times m_2$ equations. 
Therefore, we assume that $W^*=(u^*,0)\in Z_1$ and $u\to u^*$, $v\to0$ in $L^\infty(\Og)$. Recall that $u^*$ solves
\beqno{ueqnz}-\Div(\ma^{11}(u^*,0)Du^*)=\mb^{11}(u^*,0)Du^*+\mg^{11}(u^*,0)u^*.\eeq

The equation determines the $v$-stability at $W^*=(u^*,0)$ is (see \cite{dlebook1})
\beqno{per10}\barrl{\tau(-\Div(\ma^{22}(W^*)D\psi)-\mb^{22}D\psi-\Div(\ma^{21}_v(W^*)Du^*\psi) +k\psi)=}{7cm}&(\mb^{21}_v(W^*)Du^*+\mg^{22}(W^*)+k)\psi.\earr\eeq
$(u^*,0)$ is $v$-unstable if the principal eigenvalue $\tau$ of this equation is {\em greater than 1} and {\em the above equation has a positive solution $\psi$}. To be precise, $\ma^{21}_v$ should be $\tau^{-1}\ma^{21}_v$, the proof can be modified accordingly with minor changes. In this paper,  we keep this assumption for simplicity.

One has to be careful about the sizes of matrices/vectors and order of multiplications in the above equations if $m_1\ne m_2$. For simplicity one can assume $m_1=m_2$ as the calculations and assumptions below are very similar otherwise. Indeed, dropping the $x\in\Og$ and for $a=[a_{ij}]$ being a $m_1\times m_2$ matrix, $v,\psi,\fg\in\RR^{m_2}$ and $u\in\RR^{m_1}$ we adapt the following notations (which are vectors in $\RR^{m_2}$ or $m_2\times m_2$ matrix) $a_v Du \psi=a_v\psi Du =[(a_{ij})_v\psi Du_j]_{i=1}^{m_2}$, $\myprod{a_v Du \psi, D\fg}=[(a_{ij})_v\myprod{Du_j,D\fg_i}]\psi$ and  $a_vDu=[(a_{ij})_vDu_j]$.

\subsection{Technical lemmas}\label{dyntech}

Our analysis largely based on the following technical lemma on a solution  $W=(u,v)\in\RR^{m_1}\times\RR^{m_2}$ of the system \mref{syswhole} in the Introduction. The  system of $v$ is 
\beqno{per2}v_t=\Div(\ma^{21}(W)Du+\ma^{22}(W)Dv)+\mb^{21}(W)Du+\mb^{22}(W)Dv+\mg^{22}(W)v.\eeq

The $v$-unstability of a semi-trivial solution $W^*=(u^*,0)$ is determined by the assumption that the eigenvalue problem \mref{per10}
has a positive solution $\psi$ for some numbers $k$.
Inspired by this, for any $m_2\times m_2$ square matrix  $K$ we also consider the following more general system 
\beqno{per1}\barrl{\tau(-\Div(\ma^{22}(W^*)D\psi)-\mb^{22}(W^*)D\psi-\Div(\ma^{21}_v(W^*)Du^*\psi) +k\psi)=}{5cm}&(\mb^{21}_v(W^*)Du^*+\mg^{22}(W^*)+kId+K)\psi \earr
\eeq
and assume that it has a positive solution $\psi$ for some $\tau> 1$ (here, $Id$ is the $m_2\times m_2$ identity matrix). 

The following technical lemma is our starting point and it is crucial.

\blemm{Yeqn} Let $\nu\in C^1(\Og)$. Set $\ccY(t)=\iidx{\Og}{\myprod{v,\ma^{22}(W^*)\psi\nu}}$ and $\kappa=\frac{k(\tau-1)}{\tau}$. We have that $\ccY'(t)= \kappa\ccY(t)+\mC(t,u,v)+\mB(t,u,v)+\mD(t,u,v)$ with $\mC(t,u,v),\mB(t,u,v),\mD(t,u,v)$ are bounded functions given by
$$\barrl{\mC(t,u,v):= \iidx{\Og}{\myprod{\ma^{22}(W^*)D\psi,\mP_1 D\nu}}+\iidx{\Og}{\myprod{\ma^{21}_v(W^*)Du^*\psi,\mP_1D\nu}}}{.5cm}
&-\iidx{\Og}{\myprod{\ma^{21}(W)Du,\ma^{22}(W^*)D(\psi\nu)}}+\iidx{\Og}{\myprod{\mg^{22}(W)v,\ma^{22}(W^*)\psi\nu}}
\\&-\iidx{\Og}{\myprod{\ma^{22}(W)Dv,\ma^{22}_u(W^*)Du^*\psi\nu}}-\iidx{\Og}{\myprod{\ma^{21}(W)Du,\ma^{22}_u(W^*)Du^*\psi\nu}}\\
&+\iidx{\Og}{\myprod{\ma^{22}(W^*) D\psi,\mP_2\nu}}-\iidx{\Og}{\myprod{\tau^{-1}[\mg^{22}(W^*)+K]\psi,\nu\mP_1}}
\\&+\iidx{\Og}{\myprod{\ma^{21}_v(W^*)Du^*\psi,\mP_2+\ma^{22}(W)Dv}\nu}
,\earr
$$

$$\mB(t,u,v):=
-\iidx{\Og}{\nu\left(\myprod{\tau^{-1}k(1-\tau)\psi, \mP_1}-\myprod{v,\tau^{-1}k(1-\tau)\ma^{22}(W^*)\psi}\right)},
$$

$$\mD(t,u,v)=\iidx{\Og}{\myprod{\mb^{22}(W)Dv+\mb^{21}(W)Du,\psi\nu}+\myprod{\mb^{22}(W^*)D\psi-\mb^{21}_v(W^*)Du^*,\mP_1}}.$$
Here, $$\mP_1=\left[\sum_{j=1}^{m_2}\int_0^{v_j}\ma^{22}_{i,j}(u,\bar{s}_j)ds\right]_1^{m_2},\; \mP_2=\left[\sum_{l=1}^{m_1}\int_0^{v_j}(\ma^{22}_{i,j})_{u_l}(u,\bar{s}_j)dsDu_l+\sum_{l=1}^{m_2}\int_0^{v_j}(\ma^{22}_{i,j})_{v_l}(u,\bar{s}_j)dsDv_j\right]_1^{m_2}.$$

\elemm

For simplicity, {\em we will assume that $\mb\equiv0$}. The calculations below can easily be extended otherwise. In fact, the inclusion of the term $\mb$ would give rise to terms which can be handled in the same ways (see \refrem{bne01}).

\bproof Multiplying \mref{per2} with $\fg=\ma^{22}(W^*)\psi\nu$ and noting that $$D\fg=\ma^{22}(W^*)D(\psi\nu)+\ma^{22}_u(W^*)Du^*\psi\nu,$$ we integrate by parts to  have
\beqno{AAA}\barrl{\frac{d}{dt}\iidx{\Og}{\myprod{v,\ma^{22}(W^*)\psi}\nu}=-\iidx{\Og}{\myprod{\ma^{22}(W)Dv,\ma^{22}(W^*)D\psi}\nu}}{2cm}&-\iidx{\Og}{\myprod{\ma^{22}(W)Dv,\ma^{22}(W^*)\psi D\nu}\psi}-\iidx{\Og}{\myprod{\ma^{21}(W)Du,\ma^{22}(W^*)D(\psi\nu)}}\\
&-\iidx{\Og}{\myprod{\ma^{22}(W)Dv,\ma^{22}_u(W^*)Du^*\psi}\nu}-\iidx{\Og}{\myprod{\ma^{21}(W)Du,\ma^{22}_u(W^*)Du^*\psi}\nu}\\
&+\iidx{\Og}{\myprod{\mg^{22}(W)v,\ma^{22}(W^*)\psi}\nu}.
\earr\eeq

Define $$\fg =\nu\left[\sum_{j=1}^{m_2}\int_0^{v_j}\ma^{22}_{i,j}(u,\bar{s}_j)ds\right]_{i=1}^{m_2},\; \bar{s}_j=[s_i]_1^{m_2} \mbox{ with } s_i=\left\{\barr{ll}s& i=j,\\v_i&i\ne j.\earr\right.$$

Note that $$\barr{lll}D\fg&=&D\nu\left[\sum_{j=1}^{m_2}\int_0^{v_j}\ma^{22}_{i,j}(u,\bar{s}_j)ds \right]_1^{m_2}+\nu\ma^{22}(W)Dv+\\
&&\nu\left[\sum_{l=1}^{m_1}\int_0^{v_j}(\ma^{22}_{i,j})_{u_l}(u,\bar{s}_j)dsDu_l+\sum_{l=1}^{m_2}\int_0^{v_j}(\ma^{22}_{i,j})_{v_l}(u,\bar{s}_j)dsDv_j\right]_1^{m_2}\\
&=& D\nu \mP_1+\nu\mP_2+\nu\ma^{22}(W)Dv.\earr$$ 

Testing $\fg=\nu\mP_1$ with \mref{per1} which is written as
$$\barrl{-\Div(\ma^{22}(W^*)D\psi)=\mb^{22}D\psi+\Div(\ma^{21}_v(W^*)Du^*\psi) +}{7cm}&\tau^{-1}(\mb^{21}_v(W^*)Du^*+\mg^{22}(W^*)+k(1-\tau)+K)\psi\earr$$
As $\mb=0$ (for simplicity), \mref{per1} is
$$-\Div(\ma^{22}(W^*)D\psi)=\Div(\ma^{21}_v(W^*)Du^*\psi) +\tau^{-1}(\mg^{22}(W^*)+k(1-\tau)+K)\psi.$$

Denoting as in the lemma, we have 
$$\barrl{\iidx{\Og}{\myprod{\ma^{22}(W)Dv,\ma^{22}(W^*)D\psi}\nu}=-\iidx{\Og}{\myprod{\ma^{22}(W^*)D\psi,\mP_1 D\nu}}-\iidx{\Og}{\myprod{\ma^{22}(W^*)D\psi,\mP_2}\nu}}{2cm}&
\\
&-\iidx{\Og}{\myprod{\ma^{21}_v(W^*)Du^*\psi,\mP_1D\nu}}+\iidx{\Og}{\myprod{\ma^{21}_v(W^*)Du^*\psi,\mP_2+\ma^{22}(W)Dv}\nu}\\&
+\iidx{\Og}{\myprod{(\tau^{-1}\mg^{22}(W^*)\psi+\tau^{-1}k(1-\tau)\psi+K\psi),\mP_1}\nu}.
\earr$$

Therefore, combine the above two equations and  rearrange to easily get
\beqno{vCeqn}\frac{d}{dt}\iidx{\Og}{\myprod{v,\ma^{22}(W^*)\psi}\nu}=\frac{k(\tau-1)}{\tau}\iidx{\Og}{\myprod{v,\ma^{22}(W^*)\psi}\nu}+\mC(t,u,v)+\mB(t,u,v),\eeq where $\mC(t,u,v),\mB(t,u,v)$ are the quantities given in the lemma. \eproof

\brem{bne01} If $\mb\ne0$ then note that  $\mb^{21}(W^*)=0$  and the extra terms $\mD$. We also have to replace $\mg^{22}(W^*)$ by $\mb^{21}_v(W^*)Du^*+\mg^{22}(W^*)$.

\erem

{\bf A general version of \reflemm{Yeqn}:} Although we will mainly discuss the consequences of the equation for $\ccY$ in \reflemm{Yeqn} but we would like to present here a more general version of it for future investigations. 

We can replace the number $\tau,k$   in \mref{per1} by diagonal $m_2\times m_2$ constant matrices $\Theta,\LLg_*$ respectively ($\Theta,\LLg_*$ are invertible). If $(\Theta^{-1}\LLg_*-\Theta^{-1}\LLg_*\Theta)$ commutes with $\ma^{22}(W^*)$ then
the above proof can be repeated verbatim and assume that the following system has a positive solution
\beqno{Per1}\barrl{\Theta(-\Div(\ma^{22}(W^*)D\psi)-\mb^{22}D\psi-\LLg^{-1}\Div(\ma^{21}_v(W^*)Du^*\psi) +\LLg_*\psi)=}{7cm}&(\mb^{21}_v(W^*)Du^*+\mg^{22}(W^*)+\LLg_*)\psi.\earr\eeq

We then obtain a more general result.

\blemm{per1rem} Let $\psi$ be a solution of \mref{Per1} and $\ccY(t)=\iidx{\Og}{\myprod{v,\ma^{22}(W^*)\psi\nu}}$. Assume that
\bdes\item[a.1)] $(\Theta^{-1}\LLg_*-\Theta^{-1}\LLg_*\Theta)$ commutes with $\ma^{22}(W^*)$ which is symmetric.

\edes
Then we have the following equation   \beqno{per1eqna}\ccY'(t)= \iidx{\Og}{\myprod{v,(\Theta^{-1}\LLg_*\Theta-\Theta^{-1}\LLg_*)\ma^{22}(W^*)\psi}\nu}+\mC(t,u,v)+\mB(t,u,v)+\mD(t,u,v).\eeq \elemm

The condition a.1) is needed to handle the integrands in $\mB$ and the first integral on the right hand side of \mref{per1eqna}.
The term $\iidx{\Og}{\nu\left(\myprod{\tau^{-1}k(1-\tau)\psi, \mP_1}-\myprod{v,\tau^{-1}k(1-\tau)\ma^{22}(W^*)\psi}\right)}$ in $\mB$ is important. It is now 
$$\iidx{\Og}{\nu\left(\myprod{(\Theta^{-1}\LLg_*-\Theta^{-1}\LLg_*\Theta)\psi, \mP_1}-\myprod{v,(\Theta^{-1}\LLg_*-\Theta^{-1}\LLg_*\Theta)\ma^{22}(W^*)\psi}\right)}.$$

Later we will consider the situation $W\to W^*$ and thus $\mP_1\sim \ma^{22}(W^*)v$. Because of a.1),  {\em $\ma^{22}$ is symmetric} and commutes with $(\Theta^{-1}\LLg_*-\Theta^{-1}\LLg_*\Theta)$ so that it is easy to see that
$$
\iidx{\Og}{\nu\myprod{(\Theta^{-1}\LLg_*-\Theta^{-1}\LLg_*\Theta)\psi, \mP_1}}\sim \iidx{\Og}{\nu\myprod{v,(\Theta^{-1}\LLg_*-\Theta^{-1}\LLg_*\Theta)\ma^{22}(W^*)\psi}}.$$

We would also like to derive a differential inequality of the form $\ccY'\ge c_*\ccY$ for some constant $c_*>0$.

\blemm{Per2rem} In addition to a.1) we assume that

\bdes
\item[a.2)] there is $\kappa_*>0$ such that the entries $[(\Theta^{-1}\LLg_*-\Theta^{-1}\LLg_*\Theta)-\kappa_*Id]\ma^{22}(W^*)$ are nonnegative .

\edes
If $W\to W^*$ then we have $\ccY'\ge \kappa_*\ccY+\mC(t,u,v)+\mB(t,u,v)+\mD(t,u,v)$. \elemm

Concerning the first term on the right side of the equation for $\ccY$, we need that if $v,\psi>0$ for then some $\kappa_*>0$
$$\myprod{v,(\Theta^{-1}\LLg_*\Theta-\Theta^{-1}\LLg_*)\ma^{22}(W^*)\psi}\ge \kappa_*\myprod{v,\ma^{22}(W^*)\psi}.$$
This is because we assume a.2), $[(\Theta^{-1}\LLg_*-\Theta^{-1}\LLg_*\Theta)-\kappa_*Id]\ma^{22}(W^*)\psi\ge0$ for all $\psi\ge 0$ as the entries of $[(\Theta^{-1}\LLg_*-\Theta^{-1}\LLg_*\Theta)-\kappa_*Id]\ma^{22}(W^*)$ are nonnegative. 

In particular, if $\ma^{22}_{i,j}\ge0$ for all $i,j$ we can take $\Theta=\mbox{diag}[\tau_i], \LLg_*=\mbox{diag}[k_i]$ then
$$\myprod{v,\Theta^{-1}\LLg_*\Theta-\Theta^{-1}\LLg_*)\ma^{22}(W^*)\psi}=\sum_i\kappa_i v_i\ma^{22}_{i,j}(W^*)\psi_j, \quad  \kappa_i=\frac{k_i(\tau_i-1)}{\tau_i}.$$

However, in this work we will mainly consider \mref{per1}. Note that, even in this case, we still need that $\ma^{22}$ is symmetric as we see in the next section so that $\myprod{\ma^{22}v,\psi}
=\myprod{\ma^{22}\psi,v}$. But $\ma^{22}$ needs not be  symmetric in \reflemm{Yeqn}!

\subsection{Checking m.4) for \mref{syswhole} when $\tau>1$}
{\bf The Neumann boundary conditions:}  We now use a crucial fact that $\inf_{\bar{\Og}}\psi>0$, which is guaranteed by maximum principles if the homogeneous Neumann boundary condition is assumed (e.g. see \cite{Amsurvey}). Again, our crucial assumption is that the equation \mref{per1} has a positive solution $\psi$.

For simplicity, we assume first that $\mb=0$ in the statement and proof below (see \refrem{bnonzerorem} after the proof). Also, we will always assume that the $C^1$ norms of the solutions are uniformly bounded because the associated dynamical system possesses a global attractor in $C^1(\Og)$. 
We refer the readers to  \cite{dlebook1,dleJMAA,GiaS}  on how to obtain uniform bound for the derivatives $Du, Dv$ once we obtained that the solution $(u,v)$ is (uniform) H\"older continuous.

\btheo{perNthm} Assume the homogeneous Neumann boundary condition (so that $\inf_{\bar{\Og}}\psi>0$) and that $\sup_{t\ge0}\|W\|_{L^\infty(\Og)}$ and $\sup_{t\ge0}\|DW\|_{L^\infty(\Og)}$ are uniformly bounded. 

Suppose that the coefficients $\ma^{22},\ma^{21}$ satisfy the following estimates, along any positive solution $W=(u,v)$,
for some sufficiently small $\eg>0$ (in terms of $|\kappa|$, $\kappa=\frac{k(\tau-1)}{\tau}$, which is given in \reflemm{Yeqn}) 

\beqno{ma21}|\ma^{21}(W)|\le \frac{\eg v}{\sup_\Og(|Du|+|Du^*|+|D\psi|+1)},\eeq 
\beqno{ma2} |\ma^{22}(W^*)|\le \frac{\eg \inf\psi}{\sup_\Og|D\psi|+1}, \eeq
\beqno{ma3}|\ma^{21}_v(W^*)|\le \frac{\eg}{\sup_\Og|D\psi|+1},  \eeq
\beqno{ma4}|\ma^{22}_u(W^*)||D^2u^*|\le \eg,\; |\ma^{22}_u(W^*)||Du^*|\frac{|D\psi|}{\psi}\le\eg,   \eeq
\beqno{ma5} |\ma^{22}_u(W)\ma^{21}_v(W^*)||Du^*||Du|\le \eg |\ma^{22}(W^*)|,\; |\ma^{21}_v(W^*)||Du^*|\frac{|D\psi|}{\psi}\le \eg.\eeq

Moreover, suppose that $|\ma^{22}_u(u,s)|$ is small for all  $\bar{s}$ in the box $[0,v]$ ($v\in \RR^{m_2}$) and
\bdes\item[P.1)] $\ma^{22}$ is symmetric with positive entries and $\ma^{22}_{u,u}(W^*), -\ma^{21}_{u,v}(W^*)$ are positive definite matrices. Also, the off-diagonal entries of $\mg^{22}$ are nonnegative. In addition, $\ma^{22}, (\mg^{22}-G)$ are commutative ($G:=\mbox{diag}[G_i]$ with $G_i=\mg^{22}_{ii}(W^*)$).
\edes

Finally suppose that $\tau>1$, $\ma^{22}_{ij}(W^*)\ge0$ for all $i,j$, $\mg^{22}_{ij}(W^*)\ge0$ if $i\ne j$ (that is, $\ma^{22},\mg^{22}$ are {\em cooperative}) and \beqno{gper1}\tau^{-1}\sum_{i,j}|k_{ij}|+|-\tau^{-1}G_i+G_j|\le \kappa/4,\quad K=[k_{ij}] \mbox{ and } \kappa=\frac{k(\tau-1)}{\tau}.\eeq

Then there is $\eta_0>0$ such that $\limsup_{t\to \infty}(\|v(t)\|_{L^\infty(\Og)}+\|u(t)-u^*\|_{L^\infty(\Og)})\ge \eta_0$ (or equivalently $W\not\to W^*$).

\etheo

In fact, we can take $\Theta=\tau Id$ and $\LLg_*=k Id$ in a.1) of \reflemm{per1rem}, a.2) of \reflemm{Per2rem}  and obtain similar result.

For simplicity of our presentation, we assume that $\ma^{22}$ is symmetric and certain commutativity in P.1) but these assumptions can be relaxed and we refer the readers to \refrem{symmeric} and \refrem{Garem} after the proof.

\bproof The proof is by contradiction. Suppose that for any $\eta>0$, $$(\|v(t)\|_{L^\infty(\Og)}+\|u(t)-u^*\|_{L^\infty(\Og)})\le\eta \quad\mbox{ for all $t\ge0$.}$$

We now take $\nu\equiv 1$ in \reflemm{Yeqn}. The terms in $\mC$ involving $Du,Dv$ are
$$ I_1=-\iidx{\Og}{\myprod{\ma^{21}(W)Du,\ma^{22}(W^*)D\psi}},\;I_2=-\iidx{\Og}{\myprod{\ma^{22}(W)Dv,\ma^{22}_u(W^*)Du^*}\psi},$$
$$I_3=-\iidx{\Og}{\myprod{\ma^{21}(W)Du,\ma^{22}_u(W^*)Du^*\psi}},\;I_4=\iidx{\Og}{\myprod{\ma^{22}(W^*)D\psi,\mP_1}},$$
$$I_5=\iidx{\Og}{\myprod{\ma^{21}_v(W^*)Du^*\psi ,\ma^{22}(W)Dv}}, \; I_6=-\iidx{\Og}{\myprod{\ma^{21}_v(W^*)Du^*\psi,\mP_1}}.$$
Recall that $$\mP_1=\left[\sum_{j=1}^{m_2}\int_0^{v_j}\ma^{22}_{i,j}(u,\bar{s}_j)ds\right]_1^{m_2}\sim \ma{22}(W^*)v \mbox{ as $W\to W^*$}.$$

In the proof if a function $f\in L^\infty(Q)$, $Q=\Og\times[\infty)$, then with a slight abuse of notation we will simply write $\|f\|_{L^\infty(\Og\times[\infty))}$ by $\sup_Q |f|$.

Recall also that $\ccY(t)=\iidx{\Og}{\ma^{22}(W^*) v\psi}$. Here, we use the fact that $\inf_{\bar{\Og}}\psi>0$ which holds as we assume the homogeneous Neumann boundary condition here.

We are going to estimate $I_i$'s to obtain a differential inequality $\ccY'\ge \frac{\kappa}{4}\ccY$ and achieve a contradiction to complete our proof. 

It is straightforward to see that $|I_1|,|I_3|\le \eg\ccY$ (respectively $|I_4|\le \eg\ccY$) if we assume in \mref{ma21} (respectively, \mref{ma2}). 
Similarly, 
if  \mref{ma3} holds then $|I_6|\le \eg\ccY$.

Using integrations by parts for $I_2$, we obtain the integrals
\beqno{I2ints} \iidx{\Og}{\myprod{v,\ma^{22}(W)\psi\myprod{D(\ma^{22}_u(W^*)),Du^*}}}+\iidx{\Og}{\myprod{v,\ma^{22}(W^*)D(\ma^{22}(W)Du^*\psi)}}.\eeq

The first term is nonnegative because $\myprod{v,\ma^{22}(W)\psi}\ge0$ (as the involving entries are positive) and  by a.2) and $Dv^*=0$,  $\myprod{D(\ma^{22}_u(W^*)),Du^*}=\myprod{\ma^{22}_{u,u}Du^*,Du^*}\ge0$. The second term in \mref{I2ints} can be estimated by $\eg\ccY$ thanks to \mref{ma4}. We just proved that $I_2\ge \eg\ccY$.

Concerning $I_5$, we will prove by a similar argument that
$$I_5=\iidx{\Og}{\myprod{\ma^{21}_v(W^*)Du^*\psi,\ma^{22}(W)Dv}} =\iidx{\Og}{\myprod{\ma^{22}(W)^T\ma^{21}_v(W^*)Du^*\psi,Dv}}\ge \eg\ccY(t).$$

Integrating by parts, among others, yields the following integral which involves the second derivatives of the cross diffusion $\Div(\ma^{21}_v(W^*)Du^*)=\myprod{\ma^{21}_{u,v}(W^*)Du^*,Du^*}+\ma^{21}_v\Delta u^*$ 

$$\barr{lll}-\iidx{\Og}{\myprod{\ma^{22}(W)^T\Div(\ma^{21}_v(W^*)Du^*)\psi, v}}&=&-\iidx{\Og}{\myprod{\ma^{21}_{u,v}(W^*)Du^*,Du^*}\myprod{\ma^{22}(W)^T\psi, v}}\\
&&-\iidx{\Og}{\ma^{21}_{v}(W^*)D^2 u^*\myprod{\ma^{22}(W)^T\psi, v}}
\\
&&-\iidx{\Og}{\myprod{\ma^{22}(W)^T\ma^{21}_v(W^*)D^2u^*\psi, v}}.\earr$$ The first integral on the right hand side  is nonnegative because  $-\ma^{21}_{u,v}(W^*)$ is a positive definite  matrix (by a.2)) and the entries of $\ma^{22}(W)^T$ are positive {this is similar to the treatment of the first term in \mref{I2ints}}. The other integrals are less than $\eg\ccY(t)$ as $|\ma^{21}_v(W^*)||D^2u^*|$ is small by \mref{ma5} (keep in mind that $\ma^{21}(u,0)=0$).
Therefore, $I_5\ge \eg\ccY(t)$.

The assumption $|\ma^{22}_u(u,s)|\le c_*\eg$ for $s\in[0,v]$ implies that $|\mP_2|\le \eg v$ so that we can treat the integrals in $\mC$ involving $\mP_2$ and show that they can be estimate them by $\eg\ccY$.

If $\tau>1$ then $\kappa>0$. By choosing $\eg$ is sufficiently small, in terms of $\kappa$,  we see from \mref{vCeqn} that $\frac{d}{dt}\ccY(t)\ge \frac{\kappa}{2}\ccY(t)+\hat{\mC}$, where $\hat{\mC}$ consists
of terms not involving $Du$ or $Dv$ in $\mC,\mB$
\beqno{vpsiC}\barrl{\iidx{\Og}{\myprod{\tau^{-1}(\mg^{22}(W^*)+K)\psi,\mP_1 }}+\iidx{\Og}{\myprod{\mg^{22}(W)v,\ma^{22}(W^*)\psi}}+ }{4cm}&-\iidx{\Og}{\left(\myprod{\tau^{-1}k(1-\tau)\psi,\mP_1}-\myprod{v,\tau^{-1}k(1-\tau)\ma^{22}(W^*)\psi}\right)}.
\earr\eeq

This quantity  can be written as
$$\iidx{\Og}{\myprod{\psi,\left(\tau^{-1}[\mg^{22}(W^*)+K-k(\tau-1)Id]^T\mP_1-[\mg^{22}(W)-\tau^{-1}k(\tau-1)]\ma^{22}(W^*)^Tv\right)}}.$$

We note that when $W\to W^*$ (that is $u\to u^*$ and $v\to0$) $$\mP_1=\left[\sum_{j=1}^{m_2}\int_0^{v_j}\ma^{22}_{i,j}(u,\bar{s}_j)ds\right]_1^{m_2}\sim \ma^{22}(W^*)v,\;\mg^{22}(W)\to \mg^{22}(W^*).$$
Hence,  the last integral in the above \mref{vpsiC} converges to $0$ and \mref{vpsiC} tends to  the integral of

$$\myprod{\psi,\left\{\tau^{-1}[\mg^{22}(W^*)+K-k(\tau-1)Id]^T\ma^{22}(W^*)v-\ma^{22}(W^*)^T[\mg^{22}(W^*)-\tau^{-1}k(\tau-1)Id]v\right\}}$$
If $\ma^{22}$ is symmetric and commutes with $\mg^{22}$ then we can simplify it to
$$\myprod{\psi, [\tau^{-1}(\mg^{22}+K)^T-\mg^{22}]\ma^{22}v}.$$ 

Dropping $W^*$ for simplicity, we define $G:=\mbox{diag}[\mg^{22}_{ii}]$. The above integrand can be written as $G_1+G_2-\myprod{\tau^{-1}K\psi,\ma^{22}v}$ with 
$$G_1:=\myprod{\tau^{-1}G\psi,\ma^{22}v}-\myprod{Gv,\ma^{22}\psi},\; G_2:=\myprod{\tau^{-1}(\mg^{22}-G)\psi,\ma^{22}v}-\myprod{(\mg^{22}-G)v,\ma^{22}\psi}.$$ Since  $\tau>1$ and the entries of $\ma^{22}, (\mg^{22}-G)$ are nonnegative (by P.1)) and $v,\psi>0$, we have that $G_2\ge0$  so that $G_1+G_2$ is greater than $G_1$, which is

$$\tau^{-1}\sum_{i,j}\ma^{22}_{ij}\mg^{22}_{ii}v_j\psi_i-\sum_{i,j}\ma^{22}_{ij}\mg^{22}_{ii}v_i\psi_j=\sum_{i,j}(\tau^{-1}\ma^{22}_{ij}\mg^{22}_{ii}-\ma^{22}_{ji}\mg^{22}_{jj})v_j\psi_i=\sum_{i,j}(\tau^{-1}\mg^{22}_{ii}-\mg^{22}_{jj})\ma^{22}_{ij}v_j\psi_i.$$

So if $\tau^{-1}\sum_{i,j}|k_{ij}|+|-\tau^{-1}\mg^{22}_{ii}+\mg^{22}_{jj}|\le \kappa/4$ ($\kappa>0$) then we have $\frac{d}{dt}\ccY(t)\ge \frac{\kappa}{4}\ccY(t)$ (as the entries of $\ma^{22}$ are nonnegative). Because $\kappa>0$ and $\ccY(t)\ge0$, $\lim_{t\to\infty}\ccY(t)= \infty$  and this contradicts the facts that $\ccY$ is bounded on $[0,\infty)$ and $\lim_{t\to\infty}\ccY(t)= 0$.  The condition on $G$ here is assumed in \mref{gper1} of the theorem. The proof is complete. \eproof

\brem{symmeric} Note that we did not need that $\ma^{22}$ is  symmetric and there is certain commutativity in P.1). However, the entries of $\ma^{22}$ are assumed to be positive. If $\ma^{22}$ is not symmetric and does not commute, then the last integral in the above \mref{vpsiC} converges to $0$ and \mref{vpsiC} tends to  the integral of
(dropping $W^*$)
$$\myprod{\psi,\left\{[\tau^{-1}(\mg^{22}+K)^T\ma^{22}-(\ma^{22})^T\mg^{22}]v+\tau^{-1}k(\tau-1)((\ma^{22})^T-\ma^{22})v\right\}}.$$
The computation for this term would be more complicated. However, it is possible to give another set of assumptions to obtain the same assertions.

\erem

\brem{Garem} The condition P.1) is essentially to guarantee the estimate for $I_5\le \eg \ccY$ and $G_2$. It is easy to see that we can achieve the same goal if we assume that $\Div(\ma^{21}_v(W^*)Du^*)$ and the {\em off diagonal} entries of $\mg^{22}$ are small in comparison to $\kappa=k(1-\tau^{-1})$. Note also that $\ma^{21}, \ma^{21}_v$ can be small as required in \mref{ma21} but $\Div(\ma^{21}_v(W^*)Du^*)$ can be large in terms of $k$.\erem

Note the condition $\mg^{22}(W^*)+k\ge0$ later does NOT contradicts $|\mg^{22}|\le k/4$ if $\mg^{22}_{i,i}\ge 0$.

\brem{bnonzerorem} If $\mb\ne0$ then a careful checking of the calculation for $\mC$  (using \mref{per1} and \refrem{bne01}, replacing $\mg^{22}(W^*)$ by $\mb^{21}_v(W^*)Du^*+\mg^{22}(W^*)$) reveals that we have an extra term $\tau^{-1}\mb^{21}_v(W^*)Du^*\fg$ and need $$|\tau^{-1}\mb^{21}_v(W^*)Du^*|+|(1-\tau^{-1})\mg^{22}(W^*)-k\tau^{-1}(1-\tau_0)|\le \frac{\hat{\kappa}}{4}.$$
Because we can choose $\tau_0$ negatively large so that if $\hat{\kappa}$ is large then the assumptions on the terms involving $Du, Dv$ are still the same. Finally, the above condition is equivalent to
$$\frac{3k}{4}(\tau-\tau_0)+|\mb^{21}_v(W^*)Du^*|\le (\tau-1)(\mg^{22}(W^*)+k)\le \frac{5k}{4}(\tau-\tau_0)-|\mb^{21}_v(W^*)Du^*|.$$

We can allow $g^{22}(W^*)+k<0$, because we only need $\mb^{21}(W^*)Du^*+g^{22}(W^*)+k>0$ so that $\mb^{21}_v(W^*)Du^*>0$. Adding $(\tau-1)\mb^{21}_v(W^*)Du^*$ to both sides, we see that the above is
$$\frac{3k}{4}(\tau-\tau_0)+\tau\mb^{21}_v(W^*)Du^*\le (\tau-1)(\mb^{21}_v(W^*)Du^*+\mg^{22}(W^*)+k)\le \frac{5k}{4}(\tau-\tau_0)+(\tau-2)\mb^{21}_v(W^*)Du^*.$$

If $\tau>\tau_0$ then we must have $\tau>1$ again.

In addition, we have to consider the extra term
$$\mD(t,u,v)=\iidx{\Og}{\myprod{\mb^{22}(W)Dv+\mb^{21}(W)Du,\psi\nu}+\myprod{\mb^{22}(W^*)D\psi-\mb^{21}_v(W^*)Du^*,\mP_1}}.$$

Note that $\mb^{21}(W)\to \mb^{21}(W^*)=0$ as $W\to W^*$ and $Du$ is bounded, we can drop the integral of $\mb^{21}(W)Du$. The integral of $\myprod{\mb^{22}(W)Dv,\psi\nu}$ can be treated as $I_2$. Indeed, by integrating by parts we obtain the integral of $\myprod{v,D((\mb^{22})^T\psi)}=\myprod{v,(\mb^{22}_u)^TDu^*\psi+(\mb^{22})^TD\psi}$. The integral of this can be estimated by $\eg\ccY(t)$ if we assume $|\mb^{22}_u||Du^*|\le\eg,\;|\mb^{22}|(|D\psi|+1)\le\eg\inf\psi$. Similarly, the integral of 
$\myprod{\mb^{22}(W^*)D\psi,\mP_1}$ can be treated as $I_4$.

\erem

{\bf The Dirichlet boundary conditions and a local result:} The crucial fact that $\inf_{\bar{\Og}}\psi>0$ is no longer available  if homogeneous Dirichlet boundary condition is assumed. We have only that $\inf_{\bar{B}}\psi>0$ for any ball strictly contained in $\Og$.
Therefore, we can only have
\bcoro{perNcoroD} Assume the homogeneous Dirichlet boundary conditions
and other assumptions of \reftheo{perNthm}. Then for any ball $B$ strictly contained in $\Og$ there is $\eta_0(B)>0$ such that $$\limsup_{t\to \infty}(\|v(t)\|_{L^\infty(B)}+\|u(t)-u^*\|_{L^\infty(B)})\ge \eta_0(B).$$
\ecoro 

\bproof Indeed, we take $\nu$ in \reflemm{Yeqn} to be a $C^2$ function with $\nu\equiv1$ on $B$. We then have $\inf_{\bar{\Og}}\psi\nu>0$ and we can follow the proof of \reftheo{perNthm}. The proof will be the same if we replace $\psi$ by $\psi\nu$ and include the following term in the equation of $\ccY(t)=\iidx{\Og}{v\ma^{22}(W^*)\psi\nu}$
$$\iidx{\Og}{\ma^{22}(W)\ma^{22}(W^*)\myprod{Dv,D\nu}\psi}-\iidx{\Og}{\ma^{22}(W^*)\int_0^v\ma^{22}(u,s)ds\myprod{D\nu,D\psi}}.$$

As $v\to0$, the second integral goes to $0$. Applying integration by parts to the first one, we obtain
$$-\iidx{\Og}{v\Div(\ma^{22}(W)\ma^{22}(W^*)\psi D\nu)}.$$ Since $\Div(\ma^{22}(W)\ma^{22}(W^*)\psi D\nu)$ gives rise to bounded (or fixed) terms $DW,DW^*,D\psi$ and $D^2\nu$, this integral goes to $0$ too. Therefore, we still have $\ccY'\ge \frac{\kappa}{4}\ccY$ (or $\ccY'\ge \frac{\hat{\kappa}}{4}\ccY$) when $t$ is large and the proof can go on. 

\brem{vpsi} It is clear that the proof is the same if we have $v_i\psi_j\ge0$. This to say that $v,\psi$ have the same direction or $v,\psi$ are in the same cone.
\erem

\section{On the  condition $\tau>1$-the general case} \label{tausec}\eqnoset

In what follows we will take a close look at $\tau$ in \mref{per1}, which can be regarded as an eigenvalue problem for certain appropriately defined operator,  and present some conditions and  examples where the crucial assumption $\tau>1$ can be realized.

As before, let $W^*=(u^*,0)$, with $u^*$ solves \mref{ueqnz} $$-\Div(\ma^{11}(u^*,0)Du^*)=\mb^{11}(u^*,0)Du^*+\mg^{11}(u^*,0)u^*.$$

We  write  $$L\fg=-\Div(\ma^{22}(W^*)D\fg)-\mb^{22}(W^*)D\fg-\Div(\ma^{21}_v(W^*)Du^*\fg) +k\fg,$$
$$ M\fg=(\mb^{21}_v(W^*)Du^*+\mg^{22}(W^*))\fg+k\fg+K.
$$ 

In general, proving that $\tau>1$  is an eigenvalue and that its (normalized) eigenfunction is positive is a very hard problem and we can obtain an affirmative answer below based on a key assumption that $L^{-1} M$ {\em is an positive operator} and apply the famous Krein-Rutman theorem.

{\em Assume that $L^{-1} M$ is a strong positive preserving operator on an order Banach space $X$.} In what follows we will discuss some other ways to prove that $r_{L^{-1} M}>\tau_*$ for some given $\tau_*$.

\blemm{KRspectral} Assume that $L^{-1} M$ is a strong positive preserving operator on an order Banach space $X$.
If   $\tau_*\fg-L^{-1}M\fg=L^{-1}y$ for {\em some}  $\fg,y\in X$ such that $\fg>0$  and $L^{-1}y<0$, then  $r_{L^{-1} M}> \tau_*$.

Alternatively, if $L^{-1}y>0$ then $r_{L^{-1} M}< \tau_*$.
\elemm

\bproof
We recall a result from \cite[Corollary 7.27]{Z1}, with $T=L^{-1} M$ being a positive operator, that if we have $\nu x-Ty=\mu y$ and $x>0, y>0$ then $\mbox{sign}\mu=\mbox{sign}(\nu-r(T))$. Thus, if we can present just {\em a pairs of $x,y>0$} satisfying $\nu x-Ty=\mu y$ then $\mbox{sign}\mu=\mbox{sign}(\nu-r(T))$.

The equation $\tau_*\fg-L^{-1} M\fg=L^{-1} y$ can be written as $\tau_*\fg-L^{-1} M\fg=\mu(-L^{-1} y)$ with $\mu=-1$ and $(-L^{-1} y)>0$ so that $\mbox{sign}(\tau_*-r_{L^{-1} M})=\mbox{sign}(\mu)=-1$ by the above argument. Hence, $r_{L^{-1} M}> \tau_*$.

Of course, if the equation $\tau_*\fg-L^{-1} M\fg=L^{-1} y$ has a solution $\fg>0$ and $L^{-1} y>0$ then $\mbox{sign}(\tau_*-r_{L^{-1} M})=\mbox{sign}(1)=1$. In this case, $r_{L^{-1} M}< \tau_*$.  
\eproof

The alternative of the above lemma also says that and the principal eigenvalue of \mref{per1} $\tau>\tau_*^{-1}$.

\subsection{The diagonal $L$ case}
Of course, the above argument applies if we add $(kId-K)$ to $L,M$ and assume that $(L+(kId-K)^{-1}$ is positive (this is the case if $L$ is a diagonal operator as in \cite{lopez}) to show that $\tau>\tau_*$. The equation $\tau_*\fg-(L+(kId-K))^{-1}(M+(kId-K))\fg=(L+(kId-K))^{-1}y$  is equivalent to 
\beqno{te30}\barrl{-\tau_*\Div(\ma^{22}(W^*)D\fg)-\tau_*\mb^{22}(W^*)D\fg-\tau_*\Div(\ma^{21}_v(W^*)Du^*\fg)+(\tau_*-1)(kId-K)\fg=}{7cm}&\{\mb^{21}_v(W^*)Du^*+\mg^{22}(W^*)\}\fg +y.\earr\eeq
We need $\tau_*>1$ only when we want $\frac{k(\tau_*-1)}{\tau_*}$ large to obtain the conditions of \reftheo{perNthm} easier.
If $\tau_*=1$  this is
\beqno{te3}\barrl{-\Div(\ma^{22}(W^*)D\fg)-\mb^{22}(W^*)D\fg-\Div(\ma^{21}_v(W^*)Du^*\fg)=}{7cm}&\{\mb^{21}_v(W^*)Du^*+\mg^{22}(W^*)\}\fg +y.\earr\eeq

If $(L+(kId-K))^{-1}$ is positive and $y<0$ then $(L+(kId-K))^{-1}y<0$. Thus, let $\fg,y$ be such that $\fg>0$ and  $y<0$ satisfy \mref{te3} then $\tau>1$. (This is the same as saying that for some {\em positive} $y\in C(\Og)$ (or $C_0(\Og)$) and {\em negative} $\fg$ satisfying \mref{te3}). Therefore, we will achieve our goal that $\tau>1$ via \reflemm{KRspectral} if we can prove that

\bdes\item[e)] There are some {\em negative} $y\in C(\Og)$  and {\em positive} $\fg$ (in the solution space of \mref{te30}) satisfying \mref{te30} (or \mref{te3}).
\edes

We will investigate this condition in details in \refsec{diagsec}

\subsection{The non-diagonal case}\label{nondiagsec}

In general we may not have that $L^{-1},(L+kId)^{-1}$ is positive, unless $L$ is a diagonal operator as in \cite{lopez}. However, we can apply the maximum principles in \cite{MAXP} as follows.

Assume that for some $M$ we have that $L^{-1}M$ is a positive operator ($M$ may not be positive). Again, we consider the equation ($\tau_*=1$) \beqno{maineqn}\fg-L^{-1}M\fg=L^{-1}y\Leftrightarrow L\fg=M\fg+y.\eeq 

As before, in order to prove $r_{L^{-1}M}>1$ we just need to present $y,\fg$  such that $\fg>0$ and $L^{-1}y<0$ satisfy \mref{maineqn}. Since $L^{-1}y=L^{-1}MM^{-1}y$ and $L^{-1}M$ is a positive operator, we just need $\fg>0$ and $M^{-1}y<0$ to have $L^{-1}y<0$.

Thus, we consider the assumption in order to use \reflemm{KRspectral} to prove that $\tau>1$
\bdes\item[E)] There are some $y,\fg$ satisfy \mref{maineqn} such that $\fg>0$ and $M^{-1}y<0$.
\edes

Obviously, this assumption is e) if $M=Id$ and $L^{-1}$ is a positive operator.

If $M_1$ is such that $L^{-1}M_1$ is also a positive operator then
we can also try to find $\fg>0$ and $M_1^{-1}y<0$ but still obtain E).

Because $L^{-1}M_1$ is a positive operator, there are $\llg_1>0$ and $\fg>0$ (the principal eigenvalue and principal eigenfunction of $L^{-1}M_1$) such that $$L^{-1}M_1\fg =\llg_1\fg\Leftrightarrow L\fg=\frac{1}{\llg_1} M_1\fg.$$

If $\frac{1}{\llg_1}\fg<M_1^{-1}M\fg$ then we can choose $y=\frac{1}{\llg_1}M_1\fg-M\fg$ then $M_1^{-1}y=\frac{1}{\llg_1}\fg-M_1^{-1}M\fg<0$ and $L\fg=\frac{1}{\llg_1}M_1\fg=M\fg+y$. Thus E) is verified.

More generally, if we want to prove that the principal eigenvalue $r_{L^{-1}M}=\tau>\tau_*$ then we consider instead of \mref{maineqn} the equation $\tau_*L\fg=M\fg+y$. We replace $\llg_1,\fg$ in the above argument by $\frac{\llg_*}{\tau_*},\fg_*$ respectively, where $\llg_*,\fg_*$ satisfies $L^{-1}M_1\fg_*=\frac{\llg_*}{\tau_*}\fg_*$. We then assume that $\frac{\tau_*}{\llg_*}\fg_*<M_1^{-1}M\fg_*$.

In practice, we want to find $M_1$ such that the assumption that $\llg_1^{-1}\fg<M_1^{-1}M\fg$ looks as simple as possible. One of them may be that if $M_1^{-1}M$ is diagonal then, because $\fg>0$, we need only that $\llg_1^{-1}Id<M_1^{-1}M$. Of course, if we take $M_1=kM$ for some $k>0$ then $r_{L^{-1}M}=\llg_1/k$ so that the assumption $\llg_1^{-1}\fg<M_1^{-1}M\fg=k^{-1}\fg$ is exactly $r_{L^{-1}M}>1$. Therefore, a better choice is to choose $M_1=MK_d$ where $K_d$ is a nonconstant diagonal matrix.

Note that if $L^{-1}M$ is a positive operator then an easy sufficient condition for $L^{-1}M,L^{-1}M_1$ to be positive operators is that $M^{-1}M_1$ is positive (note that $M,M_1$ may not be positive). Thus, given $L,M$ such that $L^{-1}M$ is positive and we need to show that $r_{L^{-1}M}>1$ then we would try to find a simpler $M_1$ such that $M_1^{-1}M>0$ and that it is easier to find $\llg_1=\llg_{L^{-1}M_1}$ and its principal eigenfunction $\fg$ and verify that $\llg_1^{-1}\fg<M_1^{-1}M\fg$.

More generally, we let $M_1=MP$ for some $P$ such that $M_1^{-1}M=P^{-1}>0$. This is equivalent to that $P$ is monotone (i.e. $Px\ge0\Rightarrow x\ge0$). We then requires that $$L^{-1}M_1P\fg=\llg_1\fg \mbox{ and }\llg_1^{-1}\fg<P^{-1}\fg\Leftrightarrow L^{-1}M_1P\fg=\llg_1\fg \mbox{ and }P\fg>\llg_1\fg.$$

We will go back to the positivity of $L^{-1}M$ and  this condition in \refsec{newMAXP}.

\section{Diagonal $L$ case} \label{diagsec}\eqnoset

Let us recall the following well known result.

\blemm{lopez} (Lopez \cite[Theorem 3.1]{lopez}) Consider the following {\em diagonal} system for $u=[u_i]_1^m$ with homogeneous Dirichlet or Neumann boundary conditions
$$L_i(u_i) +ku_i-Ku=f_i,$$
where $L_i(\zeta)= -\Div(a_i(x)D\zeta)+b_i(x)D\zeta+c_i(x)\zeta$, a second order elliptic differential  operator with H\"older continuous coefficients, $c_0(x)+k$ is sufficiently large in terms of a given matrix $K=[k_{ij}]$ which is a cooperative matrix ($k_{ij}>0$ if $i\ne j$, we can also take $K=0$). Then a maximum principle holds. That is if $f_i>0$ for all $i$ then $u_i>0$ in $\Og$ for all $i$. 

Moreover, if  $k_{ij}>0$ if $i\ne j$, then the principal eigenvalue of $L_i(\fg_i) +k\fg_i-K\fg=\llg_1\fg_i$ is simple and has a positive  (vector valued) eigenfunction $\fg=[\fg_i]$.

\elemm

\brem{krem} Note that we can assume $c_0\ge0$ (by choosing $k$ large). The largeness of $c_0(x)+k$ also depends on the principal eigenfunctions of $L_i$'s. In fact, let $\llg_i,\psi_i$'s be the principal eigenpairs of $L_i$'s. $k$ should be sufficiently large such that (see the proof of \cite[Theorem 3.1]{lopez}) $$(\llg_i+k)\psi_i> \sum_j k_{ij}\psi_j.$$

\erem

\brem{lopezrem} The condition that $k_{ij}>0$ if $i\ne j$ is crucial for the last assertion to hold. Moreover, it is easy to see that one can replace $k$ by a diagonal matrix $\mbox{diag}[k_1,\ldots,k_m]$ with $k_i>0$ large.
\erem

In this section, we will assume that $\ma^{22}(u^*,0), \ma^{21}_v(u^*,0)Du^*$ are diagonal (so is $\Div(\ma^{21}_v(u^*,0)Du^*)$; but $\ma^{22}(u,v), \ma^{21}_v(u,v)Du$ can be non-diagonal when $v\ne0$)
and define $$L\fg=-\Div(\ma^{22}(W^*)D\fg)-(\mb^{22}(W^*)D\fg+\ma^{21}_v(W^*)Du^*)D\fg-\Div(\ma^{21}_v(W^*)Du^*)\fg,$$
$$M\fg=(\mb^{21}_v(W^*)Du^*+\mg^{22}(W^*))\fg,$$
so that $L^{-1} $ is diagonal and $(L+(kId-K))^{-1}$ is a positive operator (by \reflemm{lopez}) if $k>0$ is sufficiently large. Also by this way, we need only that the off-diagonal of $\mb^{21}_v(W^*)Du^*+\mg^{22}(W^*)$ are {\em nonnegative}.  We will provide examples for e).

In the sequel, let $$\mG=[\mg_{ij}]:=\mb^{21}_v(W^*)Du^*+\mg^{22}(W^*).$$
We will discuss different types of $\mG$.

Of course, e) will be true if there is a function $\fg>0$  such that 
\beqno{tte3a}-\Div(\ma^{22}(W^*)D\fg)-\mb^{22}(W^*)D\fg-\Div(\ma^{21}_v(W^*)Du^*\fg)<\mG\fg.\eeq

If $\fg=[\fg_i]_1^{m_1}$, then the above means
\beqno{tte3ab}-\Div(\ma^{22}(W^*)_iD\fg_i)-[\mb^{22}(W^*)]_iD\fg_i-\Div([\ma^{21}_{i,i}]_v(W^*)Du^*_i\fg_i)<
\sum_j\mg_{ij}\fg_j\quad \forall i.\eeq

\subsection{Cooperative case:}  We consider the case $\mg_{ij}>0$ if $i\ne j$. By adding $(kId-K)$ as previously described, $M+(kId-K)$ is a positive operator. We just need to show that \mref{tte3a} has a solution $\fg>0$.

Fix a constant cooperative matrix $K=[k_{ij}]$. For $k>0$ large, we now consider the first (vector valued) eigenfunction $\fg>0$ of
$$-\Div(\ma^{22}(W^*)D\fg)-\mb^{22}(W^*)D\fg-\ma^{21}_v(W^*)Du^*D\fg+(kId-K)\fg=\llg_1\fg.$$

Then \mref{tte3a} is
\beqno{tte3bcond}\llg_1\fg-\Div(\ma^{21}_v(W^*)Du^*)\fg<(\mG+kId-K)\fg.\eeq

Accordingly, if $\fg=[\fg_i]_1^{m_1}$  then the above means: for all $i$
\beqno{tt3abcond}\llg_1\fg_i-(k+\Div([\ma^{21}_{i,i}]_v(W^*)Du^*_i)\fg_i<
\sum_j(\mg_{ij}-k_{ij})\fg_j.\eeq

Our previous argument then establish that 

\btheo{coopthm} Assume that $\ma^{22},\mb^{22},\ma^{21}_v$ are diagonal at $(u^*,0)$ and $\mg_{ij}>0$ when $i\ne j$. Let $\fg=[\fg_i]$ be the positive  eigenfunction to the principle eigenvalue of 
$$-\Div(\ma^{22}(W^*)D\fg)-\mb^{22}(W^*)D\fg-\ma^{21}_v(W^*)Du^*D\fg+(kId-K)\fg=\llg_1\fg.$$

If for all $i$
\beqno{tt3abconda}\llg_1\fg_i-(k+\Div([\ma^{21}_{i,i}]_v(W^*)Du^*_i)\fg_i<
\sum_j(\mg_{ij}-k_{ij})\fg_j,\eeq
then $\tau=r_{L^{-1} M}>1$.
\etheo

The condition $\mg_{ij}>0$ when $i\ne j$ seems to be crucial because we have to use a cooperative matrix $K$ to guarantee that $L^{-1},M$ are positive (for $L^{-1}$ we just need $k$ large) and \mref{tt3abconda} will be more likely to happen when $\mg_{ij}$'s are large (because $\llg_1,\fg$ are independent of $\mg_{ij}$). However,  if $\mg_{ij}$'s are large then it is harder to establish the existence of a global attractor. If $\mg_{ij}$'s are not large enough then $\ma^{21}_v$ can have its effect depending on the size of $\fg$.

\brem{coopK} If $\mg_{ij}>0$ when $i\ne j$, we can replace \mref{tte3a} by
$$-\Div(\ma^{22}(W^*)D\fg)-\mb^{22}(W^*)D\fg-\mG\fg<\Div(\ma^{21}_v(W^*)Du^*\fg)$$
by defining $\fg>0$ differently. 

Indeed, et $k>0$ be sufficiently large.
From the last assertion of \reflemm{lopez}, we do not need the matrix $K$ and consider the first (vector valued) eigenfunction $\fg>0$ of
$$-\Div(\ma^{22}(W^*)D\fg)-\mb^{22}(W^*)D\fg-\ma^{21}_v(W^*)Du^*D\fg+k\fg-\mG\fg=\llg_1\fg.$$

The condition \mref{tte3a} is then
\beqno{tte3bcondz}\llg_1\fg<[\Div(\ma^{21}_v(W^*)Du^*)+k]\fg.\eeq

Accordingly, as $[\Div(\ma^{21}_v(W^*)Du^*)+k]$ is diagonal and $\fg>0$, the above means
\beqno{tt3abcondz}\llg_1<\Div([\ma^{21}_{i,i}]_v(W^*)Du^*_i)+k.\eeq

\erem

\brem{coopK*} If $\mg_{ij}>0$ when $i\ne j$ and for some fixed $\tau_*>1$, then we make use of the proof that $\tau>\tau^*>1$ and we can replace \mref{tte3a} by
$$-\tau_*\Div(\ma^{22}(W^*)D\fg)-\tau_*[\mb^{22}(W^*)D\fg-\ma^{21}_v(W^*)Du^*]D\fg+k(\tau_*-1)\fg-\mG\fg<\tau_*\Div(\ma^{21}_v(W^*)Du^*)\fg.$$

Let $k>0$ be sufficiently large.
From the last assertion of \reflemm{lopez}, we now consider the first (vector valued) eigenfunction $\fg>0$ of
$$-\Div(\ma^{22}(W^*)D\fg)-\mb^{22}(W^*)D\fg-\ma^{21}_v(W^*)Du^*D\fg+k(\tau_*-1)\fg-\mG\fg=\llg_1\fg.$$

The condition \mref{tte3a} is then
\beqno{tte3bcondz*}\llg_1\fg<\Div(\ma^{21}_v(W^*)Du^*)\fg.\eeq

As $\Div(\ma^{21}_v(W^*)Du^*)$ is diagonal and $\fg>0$, the above means
\beqno{tt3abcondz*}\llg_1<\Div([\ma^{21}_{i,i}]_v(W^*)Du^*_i).\eeq

It is easy to see that $\llg_1\sim k(\tau_*-1)-\min \mG_{ii}$ so that  if the diagonal entries of $\mG$ are large and its off-diagonal ones are (positive) small then $\llg_1<\kappa<k(1-\tau_*^{-1})$. Therefore, we can have $\Div([\ma^{21}_i]_v(W^*)Du^*_i)$ is small in comparison to $\kappa=k(1-\tau^{-1})
>k(1-\tau_*^{-1})$ with $k$ large. By \refrem{Garem} we see that the persistence is possible.

\erem

\bcoro{coopcoro} Assume that $\ma^{22},\mb^{22},\ma^{21}_vDu$ are diagonal at $(u^*,0)$ and $\mg_{ij}>0$ when $i\ne j$. Let $k>0$, and $\fg=[\fg_i]$ be positive eigenfunction to the principle eigenvalue of 
$$-\Div(\ma^{22}(W^*)D\fg)-\mb^{22}(W^*)D\fg-\ma^{21}_v(W^*)Du^*D\fg+k\fg-\mG\fg=\llg_1\fg.$$

If for all $i$
\beqno{tt3abcondab}\llg_1<\Div([\ma^{21}_{i,i}]_v(W^*)Du^*_i)+k,\eeq
then $\tau>1$. Of course, the above conditions could be modified accordingly to prove that $\tau>\tau_*$ for some fixed $\tau_*>1$.
\ecoro

\brem{Gcoop} One should note that in order to have $M$ positive, we need only that $\mg^{22}_{ij}\ge 0$. However, to obtain $\fg>0$, we had to use $K$ and thus assumed that $\mg^{22}_{ij}> 0$ as in the above \reftheo{coopthm} (for $L^{-1}$ we just need $k$ large). 
\erem 

Importantly, our results here seem to require that $\Div(\ma^{21}_v(W^*)Du^*)$ is large. This shoud be compare with the conditions in \reftheo{perNthm} that $|\ma^{21}_v|$ is small and that $\ma^{21}_{u,v}$ (without $Du^*$) is positive definite.

Of course, the condition\mref{tt3abcondab} in \refcoro{coopcoro} is easier to be verified than the corresponding \mref{tt3abcond} of \reftheo{coopthm} but this corollary requires stricter condition on $\mG$ (namely, $\mg_{ij}>0$ instead of $\mg_{ij}\ge0$ when $i\ne j$).

\subsection{Counterexamples} \label{counterex}

We present here some examples which show that if our conditions are violated then the persistence/coexistence woud not hold. The examples are true for the general $m_2$. First of all, let $\fg$ be a positive eigenfunction to the principal eigenvalue of the problem
$$-\Delta\fg=\llg_*\fg \mbox{ in }\Og,\quad \fg=0 \mbox{ on }\partial\Og$$
and $f$ be a function such that we have a unique positive solution $u^*$ to
$$-\Delta u=f(u)u \mbox{ in }\Og,\quad u=0 \mbox{ on }\partial\Og.$$

We consider the parabolic system $$W_t=\Div(ADW)+G(W)W \mbox{ in }\Og\times(0,\infty),\quad W=0 \mbox{ on }\partial\Og\times(0,\infty) \mbox{ and }W(x,0)=W_0(x) \mbox{ on }\Og.$$

Here, for $a=[a_{ij}], g=[c_{ij}]$ being constant $m_2\times m_2$ matrices and $m=m_2+1$ we define the $m\times m$ block  matrices
$$A=\left[\barr{cc}a&0\\0&1\earr\right],\; G=\left[\barr{cc}g&0\\0&f(u)\earr\right].$$ Of course, we assume that $A$ satisfies the usual elliptic condition.

First, assume that $a=\mbox{diag}[a_i]$, we are consider the non cross diffusion cases. We consider the functions $w_i=e^{k_i t}\fg$ and $W=\left[\barr{c}w_i\\u_*\earr\right]_i^{m_2}$ on $\Og\times(0,\infty)$. Of course, $W$ is a solution of  $W_t=\Div(ADW)+GW$ if
$k_i=-\llg_*a_{i}+\sum_j c_{ij}$
for all $i$. We see that there are choices of $a_i>0$ and $c_{ij}$ such that $k_i<0$ so that $w_i\to0$ as $t\to \infty$. The matrix $g$ can be cooperative or competitive.

Consider the case $a$ is nondiagonal, we are considering constant cross diffusion. For some numbers $k_i$'s we consider the functions $w_i=e^{k_i t}\fg$ and $W=\left[\barr{c}w_i\\u_*\earr\right]_i^{m_2}$ on $\Og\times(0,\infty)$. It is clear that $W_t=\Div(ADW)+GW$ is equivalent to
$$k_ie^{k_i t}\fg=-\llg_*\sum_j a_{ij}e^{k_j t}\fg+\sum_j c_{ij}e^{k_j t}\fg$$
for all $i$. This is true if 
$k_i=-\llg_* a_{ii}+c_{ii},\; -\llg_* a_{ij}+c_{ij}=0 \quad\forall i\ne j.$

We now choose positive $a_{ij}, c_{ij}$ such that $0<c_{ii}<\llg_* a_{ii}$ and $c_{ij}=\llg_* a_{ij}$ if $j\ne i$. Then $k_i<0$ so that $w_i\to 0$ as $t\to \infty$.

As $a_{ij}$ are constant, the eigenvalue problem in \refcoro{coopcoro} is the system
$$-a\Delta\Fg+k\Fg -G\Fg=\llg_1\Fg,$$ which has a positive solution $\Fg=[\fg]_1^{m_2}$ if $\llg_1=\llg_*\sum_j a_{ij}+k-\sum_j c_{ij}$. We can take $a_{ij}>0$ and have that $a$ is elliptic. Because  $c_{ij}=\llg_* a_{ij}>0$ if $i\ne j$ so that $g$ is cooperative. Of course,  because $\llg_*\sum_j a_{ij}-\sum_j c_{ij}=0$ and thus the condition \mref{tt3abconda} of the corollary is violated (also the diagonality) for any choice of $k$. We see that there is a solution $W$ with positive initial data $W_0$ with $v=[w_i]_1^{m_2}\to 0$. Note that persistence occurs for the corresponding reaction-diffusion system as 0 is unstable for the $v$ component. The cross diffusion may destroy this phenomenon if it is not appropriate (according to \mref{tt3abconda}). Also, we note that the cross diffusion is not diagonal at $(u^*,0)$ and therefore certain diagonality,  as we assumed on $\Div(\ma^{21}_v(W^*)Du^*)$, seems to be necessary, see also \refrem{KMP}. 

{\bf Cooperative versus competitive:}

The condition $\mg_{ij}>0$ when $i\ne j$ if $m_2>1$ seems to be crucial because we have to use a cooperative matrix $K$ to guarantee that $L^{-1},M$ are positive. 

For $m_2=1$ we do not need a matrix $K$, but a constant $k>0$, and we can consider the competitive case (see \cite{dlebook1}). This is an interesting example where $0$ is stable in the $v$ direction but by introducing appropriate (nonconstant) cross diffusions $0$ becomes unstable so that persistence can occur.

For $m_2>1$ we need a matrix $K$, however we can still consider some partly competitive cases by using a change of variables as we will discuss this matter in the next section.

This, due to our argument, may reveal a natural phenomenon. Even if $u$ is harmful to one individual $v$ ($m_1=1$), $v$ can adjust its movement ($\Div(\ma^{21}_vDu^*)$) accordingly so that it cannot be driven to extinction. However, if $v$ belong to a group of species ($m_2>1$) which is trying to destroy each others then a component $v$ cannot control others by its movement to survive!

\subsection{Partial competitive case:}  We consider the case $\mg_{ij}\le0$ for some $i\ne j$. 
By appropriate change of variables, we will reduce this case to the cooperative ones considered previously.

In particular, we assume that $\mG$ is a {\em block} matrix
$$\mG(u^*,0)=\left[ \begin {array}{cc} A&B\\ \noalign{\medskip}C&D\end {array}
\right]
$$
where $A,D$ are square positive matrices of sizes $k,l$ and the off-diagonal entries of $B,C$ are nonnegative. We define \beqno{Pmatrix}\left[ \begin {array}{cc} Id_k&0\\ \noalign{\medskip}0&-Id_l\end {array}
\right].\eeq

It is natural to assume that the components partiticipating in the process will react to each others in a tit for tat way. That is, the variables in $v\in\RR^{m_2}$ can be divided into two competing goups but they support each others in theirs owns. Thus, the symmetric (across the diagonal of $\mG$) entries will have the same sign. Using permutation matrices, we can always assume $\mG$ to have this form.

We make use of a change of variables $\bar{v}=Pv$ then the system for $v$
becomes $$v_t=\Div(\bar{\ma}^{21}Du+\bar{\ma}^{22} Dv)+\bar{\mb}^{21} Du+\bar{\mb}^{22} Dv+\bar{\mG} v$$
where $\bar{\ma}^{ij}(u,v)=P^{-1}\ma^{ij}(u,Pv)P$, 
$\bar{\mG}(u,v)=P^{-1}\mG(u,Pv)P$
and $\bar{\mb}$ is defined similarly. We will see that the previous argument for the cooperative can apply here.  

The new system still satisfies the normal ellipticity. Furthermore, if $\ma^{22}(u^*,0), \ma^{21}_v(u^*,0), P$ are diagonal then the matrices $$\bar{\ma}^{22}(u^*,0)=P^{-1}\ma^{22}(u^*,P0)P=\ma^{22}(u^*,0) \mbox{ and } \bar{\ma}^{21}(u^*,0)=(P^{-1}\ma^{21}(u^*,P0)P)_v=P\ma^{21}_v(u^*,0)$$ as $P^{-1}=P$, $P^2=Id$ and
that $\bar{\ma}^{21}_v(u,v)=P\frac{\partial}{\partial Pv}\ma^{21}(u,Pv)P^2=P\frac{\partial}{\partial Pv}\ma^{21}(u,Pv)$ are diagonal. 
Similarly, $\bar{\mb}^{21}_v(u^*,0)=P\mb^{21}_{v}(u^*,0)$. Therefore, $\bar{\ma}^{21}_v(u^*,0),\bar{\mb}^{21}_v(u^*,0)$ are also diagonal. It is easy to see that
$$\bar{\mG}=P^{-1}\mG P=[\bar{\mg}_{ij}]=\left[ \begin {array}{cc} A&-B\\ \noalign{\medskip}-C&D\end {array}
\right].$$

Thus, the off-diagonal entries of $P^{-1}\mG(u^*,P0)P$  are nonnegative because those of $A,D,-B,-C$ are. Importantly, the above matrix (and $\mG$) does not have to be symmetric for P.1) to be true. 

We then have the following consequence of \reftheo{coopthm} and leave the corresponding statement of \refcoro{coopcoro} to the readers.
\btheo{compthm} Assume that $\ma^{22},\mb^{22},\ma^{21}_v$ are diagonal at $(u^*,0)$ and the block matrix $$ \mG(u^*,0)=  \left[ \begin {array}{cc} A&B\\ \noalign{\medskip}C&D\end {array}
\right]
$$
with $A,D$ are square positive matrices of sizes $k,l$ and the off-diagonal entries of $B,C$ are negative.

Let $\fg=[\fg_i]$ be positive eigenfunction to the principle eigenvalue of ($P$ is defined in \mref{Pmatrix})
$$-\Div(\ma^{22}(W^*)D\fg)-P\mb^{22}(W^*)D\fg-P\ma^{21}_v(W^*)Du^*D\fg=\bar{\llg}_1\fg.$$

If for all $i$
\beqno{tt3abcondbc}\bar{\llg}_1\fg_i-\Div([P\ma^{21}_i]_v(W^*)Du^*_i)\fg_i<
\sum_j\bar{\mg}_{ij}\fg_j,\eeq
then $\tau>1$.
\etheo

In block form, let $$P=\twomatrix{P_1&0\\0&P_2},\; \mg=\twomatrix{A&B\\C&D}.$$ Then $$P^{-1}\mg P=\twomatrix{P_1^{-1}AP_1 &P_1^{-1}BP_2\\ P_2^{-1}CP_1 & P_2^{-1}DP_2}.$$

The completely competitive or Prey-predator cases are harder as we discussed in \cite[section 3]{MAXP}.

\section{Non diagonal case}\label{nondiag}\eqnoset

We now continue with the argument in \refsec{nondiagsec} to prove that $\tau>1$. Firstly, one of the most crucial assumptions there is that $L^{-1}M$ is a strongly positive operator. We will check this assumption here for cross diffusion systems. As we have seen that for diagonal systems, the reaction parts could not be completely competitive. However, if we introduce cross diffusions then we can allow the reaction to be competitive but still obtain that $\tau>1$. Examples will be presented to show the contrast effect.

\subsection{$L^{-1}\mM$ is strongly positive - the nondiagonal case}\label{newMAXP}

\newc{\ccL}{{\cal L}}
\newc{\hccL}{\hat{{\cal L}}}
\newcommand{\Pcoop}{\mathbf{P_{coop}}}
\newc{\mbg}{\mathbf{m_{bg}}}
\newc{\ccM}{{\cal M}}

We make use of new max principles to establish first that $L^{-1}\mM$ is strongly positive. We recall the following result in \cite{MAXP} (in order to make the notation less confusing we change the notation $K$ in that paper to $\bar{K}$ and introduce the following class of matrices ${\cal M}$). 

\bdefi{Def} For a given collection of constant vectors $\hat{k}=\{\hat{k}_1,\ldots,\hat{k}_n\}$, $\hat{k}_i\in\RR^i$,  we define ${\cal M}(\hat{k})$ to be the collection of matrices $\mBB$ such that, in block form, their main diagonal blocks are of the form
$$\mat{\ag_i&\bg_i\\\cg_i&\dg_i}\mbox{ such that $\bg_i=\ag_i \hat{k}_i$}$$ for some invertible matrices $\ag_i$, $\hat{k}_i\in\hat{k}$, row vectors $\cg_i$ and constants $\dg_i$.

Also, we denote by $\MLT{n}$ the set of lower triangular matrices of size $n\times n$.
\edefi

We proved in \cite[Theorem 6.1]{MAXP} that (with a slight abuse of notation in \mref{blockgencond} and the sequel,  $c\hat{k}\bar{k}$ means $\myprod{c,\hat{k}}\bar{k}$)

\btheo{GenMPMatTKRnew} Assume that $\mA,\mB,\mT$ are square matrices and their main blocks matrices  satisfy 
\beqno{blockgencond}\left\{\barr{l} \mat{a&b\\c&d}\mbox{ satisfy $-a\hat{k}+b=-c\hat{k}\bar{k}+d\bar{k}$, where $\hat{k},\bar{k}$ are some constant vectors,}
\\
\mT=\mat{\ag&\bg\\\cg&\dg} \mbox{ with $\bg=\ag\hat{k}$ for the same vectors $\hat{k}$.}
\earr\right.\eeq
That is, $\mT\in {\cal M}(\hat{k})$. In addition, assume further that $a-\frac{1}{\dg}b\cg-\bar{k}(c-\frac{1}{\dg}d\cg)$ is invertible.

Then there is a  matrix $\mBB$ such that $\mBB\mA\mT^{-1},\mBB\mB\mT^{-1}\in \MLT{n}$.

Furthermore,  assume that $\mBB\mA\mT^{-1}=\mL^{-1}\mA_d$, $\mA_d=\mbox{diag}[\hat{a}_{1,1}\cdots,\hat{a}_{n,n}]$ with the functions $\hat{a}_{i,i}>0$ for some $\mL\in \MLT{n}$.

Suppose that  the matrix  $\hat{\mC}=D(\mL\mBB)\mA D\mT^{-1}+\mL\mBB\mB D\mT^{-1}$ is {\em cooperative}.

Denote $\mC=-(D(\mL\mBB)\mA\mT^{-1}+\mL\mBB\mB\mT^{-1})=[\mathbf{c}_{i,j}]$. Suppose that $\mC$ is diagonal (we can prove that $\mC$ is lower triangular with such choices of $\mL, \mBB,\mT$, see \cite{MAXP} for details)

Consider the inverse map $\hccL^{-1}$, which exists if $k$ is sufficiently large, of the map associated to the system 
\beqno{LMsys1}\left\{\barr{ll} -\Div(\mA DW) +\mB DW+ kW=F &\mbox{in $\Og$,}\\ W=0&\mbox{on $\partial\Og$}.\earr
\right.\eeq   Suppose that there is some positive matrix $\mP_{pos}$ and cooperative matrix $\Pcoop$ such that  for $\kappa$ sufficiently large
\beqno{MPpos} \kappa\mT^{-1}\mP_{pos}^{-1}\mT> \mT^{-1}\mP_{pos}^{-1}(\Pcoop\mT+k(\mL\mBB)).\eeq

Define the matrices $$ 
\mM:=(\mL\mBB)^{-1}[\mP_{pos}+\kappa Id -\Pcoop]\mT-k Id.$$ Then, $\hccL^{-1}\mM$ is a {\em strongly positive operator}.

\etheo

As we remark in \cite{MAXP}, if $\mA$ satisfies \mref{blockgencond} then we can choose $\mT$ to be a constant matrix without loss of generality but the choice of $\mM$ may be a bit restrictive.

We can summarize the main ideas of the proof in \cite{MAXP} as follows. We make a change of variable $v=\mT W$ and then multiply the system for $v$ with the matrix $\mBB$ and see that we can first convert it into a lower triangular system if and only if $\mBB\in {\cal M}({k})$ for some constant vectors $k$'s and $\mBB$ satisfies the properties in the theorem. Because the system is now lower triangular, we can find $\mL\in \MLT{n}$ such that when we multiply the resulted system with $\mL$ then we obtain a diagonal system for which \reflemm{lopez} can apply. The calculation is long but the proof is straightforward. Importantly, we proved that it is sufficient and necessary that $\mA,\mB$  must satisfy \mref{blockgencond} so that the described process can be carried out. Also, a mere change of variable (using $\mT$) is not sufficient to convert the system to a diagonal one.

We present here some remarks which show that such matrices can be found such that \reftheo{GenMPMatTKRnew} can apply in some nontrivial cases.

The following simple result (\cite[Lemma 6.2]{MAXP}) presents some easy cases concerning $\mC,\hat{\mC}$.

\blemm{LBTspecial} $\mC$ is diagonal in the following cases
\bdes \item[i)] $\mL\mBB$ is a constant matrix and $\mB=0$.
\item[ii)] $\mT$ is a constant matrix and $\mL\mBB$ is diagonal.
\item[iii)] $\mL\mBB, D\mT\mT^{-1}$ are diagonal.
\edes

In addition, $\hat{\mC}=0$ in cases i), ii) and $\hat{\mC}$ is diagonal in case iii).
\elemm

\brem{DTTdiag} We can completely describe $\mT$ such that $D\mT\mT^{-1}$ is diagonal. We can prove that $\mT$ must be of the form $\mbox{diag}[\dg_1,\ldots,\dg_n]C$, where $\dg_i$'s are functions and $C$ is a constant matrix (\cite[Remark 5.6]{MAXP}). We see that ii) is a special case of iii).
On the other hand, since $\mL$ is lower triangular and $\mL\mBB$ is diagonal we must have that $\mBB$ is lower triangular. In fact, \reftheo{GenMPMatTKRnew} shows that $\mBB$ must belong to ${\cal M}(k)$, see \refdef{Def},  for some constant vectors $k$ (see \cite[Lemma 5.4]{MAXP}). The main blocks of $\mBB$ must be of the forms  $$\mat{\ag&\bg\\\cg&\dg} \mbox{ $\bg=\ag k$ for some constant vectors $k$}$$ and its is easy to show that $\mBB$ is lower triangular if and only if $k$'s are zero vectors. Thus, in case iii), $\mA,\mB$ must be of the forms $\mD_1C_1$, $\mD_2C_2$ where $\mD_i, C_i$ are respectively diagonal and constant ({\em full}) matrices. Note that we can first determine $\mT=\mD C$, $\mL\mBB, \mA_d$  from $\mD_1,C_1$ so that the entries of the diagonal matrices $\mA_d$  are positive. Thus, \reftheo{GenMPMatTKRnew} can apply for such $\mA,\mB$. 
\erem

It is easy (see \cite{MAXP}) to show that for any given constant vectors $\hat{k}$'s we can construct $\mA$, then $\mB$, from matrices $\mBB,\mT\in {\cal M}(\hat{k})$, $\mL\in\MLT{n}$ and diagonal ones $\mA_d,\mB_d$.

We now try to verify E) of \refsec{nondiagsec} by showing that there is $\fg>0$ (in the solution space) such that for some $\tau>1$ \mref{per1} holds, this is the most important starting point of \refsec{dynevol}. 

Obviously, from the choices of $\mL,\mBB,\mT$ we can also write \mref{MPpos} as
\beqno{MPposa} \kappa\mT^{-1}\mP_{pos}^{-1}\mT> \mT^{-1}\mP_{pos}^{-1}(\Pcoop\mT+k\mA_d\mT\mA^{-1}).\eeq
On the other hand $\mM$ depends mildly on $\mA,\mB$ because besides $\mL\mBB,\mT$ we are free to choose $\mA_d,\mB_d$. Also, the conditions of \reftheo{GenMPMatTKRnew} seem to be technical but we will prove in \refcoro{LBTspecCOND1} they can be all verified by special choices of $\mL\mBB,\mT$ as in \reflemm{LBTspecial}.

\brem{mTPos} If we can write $\mT=\mP_{pos}^{-1}\hat{\mP}_{pos}^{-1}$ for some positive matrix $\hat{\mP}_{pos}$ then it is easy to see that the second inequality in \mref{MPpos} is $(\kappa Id-\Pcoop)\mP_{pos}^{-1}\hat{\mP}_{pos}^{-1}>k\mL\mBB$.
\erem

The equation \mref{per1} in \refsec{dynevol} is,
$$\barrl{\tau(-\Div(\ma^{22}(W^*)D\psi)-\mb^{22}(W^*)D\psi-\Div(\ma^{21}_v(W^*)Du^*\psi) +k\psi)=}{6cm}
&(\mb^{21}_v(W^*)Du^*+\mg^{22}(W^*)+kId+K)\psi.\earr
$$
For simplicity, we will assume that $\ma^{21}_v\equiv0$ in the sequel. The case $\ma^{21}_v\ne0$ is easy as we will see in \refrem{a21rem} later.

Let $\mA=\ma^{22}(W^*)$, $\mB=-\mb^{22}(W^*)$ and $\mP_{pos}, \Pcoop$ be some matrices such that  $\mP_{pos}$ is positive;  $\Pcoop$ is cooperative. Assume that  for $\kappa$ sufficiently large  \mref{MPpos} holds.

Assume that $\mA,\mB,\mT$ are square matrices and their main blocks matrices  satisfy 
\mref{blockgencond} described in \reftheo{GenMPMatTKRnew}. Based on $\mA$ we can find the matrices $\mBB,\hccL$ as in \reftheo{GenMPMatTKRnew} and define accordingly
$\mM:=(\mL\mBB)^{-1}[\mP_{pos}+\kappa Id -\Pcoop]\mT-k Id$.

Assume that  we can write
$\mb^{21}_v(W^*)Du^*+\mg^{22}(W^*)+kId+K=\mM$ then \mref{per1} can be rewritten as  $\tau\hccL \psi =\mM\psi$ with $\hccL,\mM$ described in \reftheo{GenMPMatTKRnew}: 
$ \tau(-\Div(\mA D\psi) +\mB D\psi+ k\psi)=\mM\psi$.
For the sake of brevity we also denote 
$$\mbg:=\mb^{21}_v(W^*)Du^*+\mg^{22}(W^*)+K$$ so that $\mM=\mbg+kId$.
By combining with the above definition of $\mM$, we must have $$\mbg=(\mL\mBB)^{-1}\mP_{pos}\mT+\kappa(\mL\mBB)^{-1}\mT-(\mL\mBB)^{-1}\Pcoop\mT-kId,$$ which is \beqno{mbgP}\mbg:=\mb^{21}_v(W^*)Du^*+\mg^{22}(W^*)+K=(\mL\mBB)^{-1}[\mP_{pos}+\kappa Id-\Pcoop]\mT-kId.\eeq

In order to apply the theorem to prove that $\hccL^{-1}\mM$ is strongly positive (so that $\tau>0$), we need to check the key condition  of \reftheo{GenMPMatTKRnew} on $\hat{\mC}$ and assume \mref{MPpos}. That is, 
\beqno{COND1} \hat{\mC}\mbox{ is {\em cooperative} and }\kappa\mT^{-1}\mP_{pos}^{-1}\mT> \mT^{-1}\mP_{pos}^{-1}(\Pcoop\mT+k(\mL\mBB)).
\eeq

\brem{a21rem}
If $\ma^{21}_v(W^*)\ne0$, then we redefine  $\mB=-\mb^{22}(W^*)-\ma^{21}_v(W^*)Du^*$, replace $K$ by $K+\Div(\ma_v^{21}(W^*)Du^*)$ and use $\mA$ as before.

\erem

We now use present some special cases where all conditions of \reftheo{GenMPMatTKRnew} are verified so that one can apply it to show that the key eigenvalue problem \mref{per1} in \refsec{dynevol} has a solution $\psi>0$ for some $\tau>0$.

In particular, we use consider the cases described in \reflemm{LBTspecial} to have the assumptions on $\mC, \hat{\mC}$ of \reftheo{GenMPMatTKRnew} fulfilled. The only matter left is prove that $\Pcoop, \mP_{pos}$ exist and satisfy \mref{COND1}. We have the following result in \cite{MAXP}.

\blemm{LBTspecCOND1lem} (\cite[Corrollary 6.4]{MAXP}) Assume the cases i)-iii) of \reflemm{LBTspecial}. There are a positive matrix $\mP_{pos}$ and a cooperative matrix $\Pcoop$ such that all conditions of \reftheo{GenMPMatTKRnew} are verified so that $\hccL^{-1}\mM$ is a {\em strongly positive operator}.

\elemm

{\em We conclude that \mref{per1} has a solution $\psi>0$ for some $\tau>0$. The above results  partially verify the key assumption of \refsec{nondiagsec} that $\hccL^{-1}\mM$ is strongly positive (we still have to show that $\tau>1$  as required).}

\brem{MLBTrem} From the remarks in \cite{MAXP}, we can choose $\mL\mBB,\mT>0$ appropritately and verify the most crucial assumption \mref{COND1} of \reftheo{GenMPMatTKRnew} $\kappa\mT^{-1}\mP_{pos}^{-1}\mT> \mT^{-1}\mP_{pos}^{-1}(\Pcoop\mT+k(\mL\mBB))$.

On the other hand, if $\mM=(\mL\mBB)^{-1}(\mP_{pos}+\kappa Id-\Pcoop)\mT-kId$ is completely competitive then we can choose $\mL\mBB>0$ and $-\mT>0$ such that \mref{COND1} holds
and that the cases of \reflemm{LBTspecial} and \reflemm{LBTspecCOND1lem} can be used here so that all conditions of \reftheo{GenMPMatTKRnew} are satisfied.

\erem

We also have the following result in \cite{MAXP} which shows that $\hccL$ exists if $\mM$ is properly given.
\blemm{Mgiven}
Let $\nu_*=\pm1$ and $\mL\mBB$ be satifying one of the cases of \reflemm{LBTspecial} with $\mL\mBB>0$ or its diagonal entries are positive. Suppose that $$\mM=(\mL\mBB)^{-1}\mP_{pos}\mT+(\mL\mBB)^{-1}(\kappa Id-\Pcoop)\mT-kId=(\mL\mBB)^{-1}\mP_{pos}\mT+\nu_*\mM_*,$$ where $\mM_*=\nu_*[(\mL\mBB)^{-1}(\kappa Id-\Pcoop)\mT-kId]$ for some $k,\kappa>0$ and matrices $\mP_{pos}>0$, cooperative $\Pcoop$, and $\mT$  such that $\nu_*\mT>0$ and $D\mT\mT^{-1}\in \MLT{n}$.

If $\mM_*>0$ then there is $\hccL$ as in \reftheo{GenMPMatTKRnew} such that $\hccL^{-1}\mM$ is strongly positive. In addition, we have that $(\mL\mBB)^{-1}\mP_{pos}\mT$
is positive if $\nu_*=1$ and negative if $\nu_*=-1$.

\elemm

Because $D\mT\mT^{-1}$ is diagonal,  $\mT=\mbox{diag}[\dg_1,\ldots,\dg_n]C$, where $\dg_i$'s are functions and $C$ is a constant matrix (see \refrem{DTTdiag}). Consider the simple case $\dg_i=\dg$ and the diagonal entries $\mu_i$'s of $\mL\mBB$ are positive (the case $C$ is diagonal is trivial as $\hccL$ will be diagonal too and we do not have cross diffusion). Then, in order $\mM=(\mL\mBB)^{-1}\mP_{pos}\mT+(\mL\mBB)^{-1}(\kappa Id-\Pcoop)\mT-kId$ we must have  some $\mP_{pos},\Pcoop$ such that
$\mM-kId=\mbox{diag}[\ag_i](\mP_{pos}+\kappa Id-\Pcoop)C$ where $\ag_i=\dg\mu_i^{-1}$ which have the same sign.

In the simplest case, if we choose $\ag_i\sim 1$ and $C\sim Id$ then $\mM\sim (\kappa-k)Id+\mP_{pos}-\Pcoop$ and $\mM_*\sim \nu_*[(\kappa-k)Id-\Pcoop]$. We then see that if $\mM\sim \kappa_*Id+M$ for some $\kappa_*>0$ (we can choose $\kappa>k>0$) and $\mP_{pos},\Pcoop$ are the positive and negative parts of $M$ with $\kappa_*Id-\Pcoop$ is negative then we choose $\nu_*=-1$ and $\hccL$ exists. In this case, $\mM$ can be completely competitive and a maximum principles is available if we have certain appropriate cross diffusion. Meanwhile, it is easy to construct an example where this is not the case when the system is diagonal.

\subsection{$\tau>1$ and the effects of cross diffusion} \label{crossdiffexample}

We can use \reftheo{GenMPMatTKRnew}, and various choices of $M,M_1$ to show that $\tau>0$. 
Let us recall the settings we assumed in order to apply \reftheo{GenMPMatTKRnew} to \mref{per1}. We assume first that  \beqno{COND1a} \hat{\mC}\mbox{ is {\em cooperative} and }\kappa\mT^{-1}\mP_{pos}^{-1}\mT> \mT^{-1}\mP_{pos}^{-1}(\Pcoop\mT+k(\mL\mBB))
\eeq and then set $\mM=\mbg+kId$ with $\mbg:=\mb^{21}_v(W^*)Du^*+\mg^{22}(W^*)+K$ to obtain  $\hccL^{-1}\mM$. We need $\hccL^{-1}\mM$, the operator
in \reftheo{GenMPMatTKRnew}, is strongly positive to show first that $\tau>0$. Thus, \beqno{mbgPa}\mbg=(\mL\mBB)^{-1}[\mP_{pos}+\kappa Id-\Pcoop]\mT-2kId.\eeq

We then apply the argument in \refsec{nondiagsec} to check the condition E) for $L=\hccL$  to show that $\tau>1$ and \mref{per1} holds. These examples show that the introduction of cross diffusion does induce persistence which may not occur in diagonal cases. One should also note that the argument in \refsec{dynevol} do not require that $K$ is a constant matrix.

{\bf When ${\cal L}$ is given:} It is easy to construct the reactions such that $\tau>1$. In this case $\mL,\mBB,\mT$ and $\mA_d$ are given. Let $\Pcoop$ (resp. $\mP_{pos}$) be a cooperative (resp. positive) matrix  satisfying \mref{COND1a} (see also \reflemm{LBTspecCOND1lem}).

Let  $\mM_1= (\mL\mBB)^{-1}[\mP_{pos}+\kappa Id-\Pcoop]\mT-kId$. By \reftheo{GenMPMatTKRnew} $\hccL^{-1}\mM_1$ is strongly positive  so that there is $\llg_*,\fg_*>0$ such that 
$$\hccL^{-1}\mM_1\fg_*=\llg_*^{-1}\fg_* \Leftrightarrow \mM_1\fg_*=\llg_*^{-1}\hccL\fg_*.$$

Take $\mM, \mb^{21}_v(W^*)Du^*$ and $\mg^{22}$ be such that $\mM\fg_*=\llg^{-1}\llg_* \mM_1\fg_*$. There are many such $\mb^{21}_v, \mg^{22}$. That is, $$(\mb^{21}_v(W^*)Du^*+\mg^{22}+K+kId)\fg_*= \llg^{-1}\llg_*(\mL\mBB)^{-1}[\mP_{pos}+\kappa Id-\Pcoop]\mT\fg_*-k\fg_*.$$  Then $\mM_1\fg_*=\llg_*^{-1}{\cal L}\fg_*\Rightarrow\mM\fg_*=\llg^{-1}{\cal L}\fg_* $.  If $\llg<1$ then \mref{per1} have a solution $\psi=\fg_*>0$ and $\tau=\llg^{-1}>1$.

The converse problem is harder. That is, if the reaction is given then we can design the diffusion in a way that $\tau>1$. If the reaction is of certain form then we can use \reflemm{Mgiven} and usual scalings to get the desired. On the other hand, we can scale the domain too. We have the following {\em partial} answer.

{\bf When $\mb^{22}=0$, and $\ma^{21}, \mb^{21},\mg^{22}$ are given on $\Og$:} We need to find $\mL,\mBB,\mT$ and $\mA_d$ to determine $\mA$ (then  $\mB=-\ma^{21}_v(W^*)Du^*$ has the structure determined by $\ma^{22}$ as in the theorem to give $\hccL$) such that \mref{COND1a} holds first. Of course, \reflemm{LBTspecCOND1lem} can be used here to obtain $\hccL,\mM$. 

\newc{\mmm}{\mathbf{m}}

Note that for any matrix $\mmm$ if we  choose $K=(\mL\mBB)^{-1}[\mP_{pos}+\kappa Id-\Pcoop]\mT-2kId-\mmm$ then there are $\fg_*>0, \tau>0$ independent of $\mmm$ such that $\hccL^{-1}\mM\fg_*=\tau\fg_*$. That is
\beqno{keytau}\tau[-\Div(\ma^{22}(x)D\fg_*(x))-\ma^{21}_v(x)Du^*(x)D\fg_*(x)]=[\mmm+K+kId]\fg_*(x).\eeq

For such fixed $\fg_*,\tau$  above we can choose $R>0$ such that $R^{-2}\tau>1$. Define $\hat{x}=Rx$ and write $\hat{\ag}(x)=\ag(Rx)$ for any function matrix $\ag$. Let $\Og^R=\{x\,:\, Rx\in\Og\}$.

By elementary calculus we have
$\tau[-\Div(\hat{\ma}^{22}(x)D\hat{\fg}_*(x))-\hat{\ma}^{21}_v(x)D\hat{u}^*(x)D\hat{\fg}_*(x)]$ is $$\tau R^2[-\Div_{\hat{x}}(\ma^{22}(\hat{x})D_{\hat{x}}\fg_*(\hat{x}))-\ma^{21}_v(\hat{x})D_{\hat{x}}u^*(\hat{x})D_{\hat{x}}\fg_*(\hat{x})]$$

So that from \mref{keytau} with $x$ being $\hat{x}$ and $\mmm(\hat{x})=\Div_{\hat{x}}(R^2\ma^{21}_v(\hat{x})D_{\hat{x}}u^*(\hat{x}))+R\mb^{21}_v(\hat{x})D_{\hat{x}}u^*(\hat{x})+\mg^{22}(\hat{x})$, we have by chain rules and combining the above calculations that
$$\tau[-\Div(\hat{\ma}^{22}(x)D\hat{\fg}_*(x))-\hat{\ma}^{21}_v(x)D\hat{u}^*(x)D\hat{\fg}_*(x)]=R^2[\mmm(\hat{x})+K+kId]\fg_*(\hat{x}).$$

Thus,
$$R^{-2}\tau[-\Div(\hat{\ma}^{22}(x)D\hat{\fg}_*(x))-\hat{\ma}^{21}_v(x)D\hat{u}^*(x)D\hat{\fg}_*(x)]=[\mmm(x)+K+kId]\hat{\fg}_*(x).$$
where $\mmm(x)=\Div(\ma^{21}_v(Rx)Du^*(Rx))+\mb^{21}_v(Rx)Du^*(Rx)+\mg^{22}(Rx)$. Therefore, \mref{per1} holds with $\tau$ being $R^{-2}\tau>1$ and $\psi(x):=\hat{\fg}_*(x)=\fg_*(Rx)>0$ for $x\in \Og^R$. 

Thus, for any given $\ma^{21}_v,\mb^{21}_v,\mg^{22}$ one has to combine an appropriate cross diffusion and scaling to have persistence on a suitably domain scaled from the original one.

{\bf A simple example when $\mA$ is a constant matrix:} We recall the counter example in \refsec{counterex} where we consider  $\mA=\mbox{diag}[a_i]$, $M=[c_{ij}]$ and  the functions $w_i=e^{k_i t}\fg$ and $W=\left[w_i\right]_i^{m_2}$ on $\Og\times(0,\infty)$. Of course, $W$ is a solution of  $W_t=\Div(\mA DW)+MW$ if
$k_i=-\llg_*a_{i}+\sum_j c_{ij}$
for all $i$. We see that there are choices of $a_i>0$ and $c_{ij}$ such that $k_i<0$ so that $w_i\to0$ as $t\to \infty$. The matrix $g$ can be cooperative or competitive. Thus, we cannot have $\tau>1$ in this case.

For some numbers $k_i$'s we consider the functions $w_i=e^{k_i t}\fg$ and $W=\left[\barr{c}w_i\earr\right]_i^{m_2}$ on $\Og\times(0,\infty)$. It is clear that $W_t=\Div(ADW)+GW$ is equivalent to
$$k_ie^{k_i t}\fg=-\llg_*\sum_j a_{ij}e^{k_j t}\fg+\sum_j c_{ij}e^{k_j t}\fg$$
for all $i$. This is true if 
$k_i=-\llg_* a_{ii}+c_{ii},\; -\llg_* a_{ij}+c_{ij}=0 \quad\forall i\ne j.$
We now choose positive $a_{ij}, c_{ij}$ such that $0<c_{ii}<\llg_* a_{ii}$ and $c_{ij}=\llg_* a_{ij}$ if $j\ne i$. Then $k_i<0$ so that $w_i\to 0$ as $t\to \infty$. This has been discussed in our previous counterexamples when  the condition \mref{tt3abconda} of the  \refcoro{coopcoro} is violated (also the diagonality). Note that this situation happen even when $\llg_*$ is large by scaling the domain $\Og$ (we can choose $a_{ij}=c_{ij}/\llg_*$ very small if $i\ne j$ and $a_{ii}$ large for given $c_{ij}$'s).

Thus, we can have persistence if $\mA$ is full but appropriately designed but no persistence when $\mA$ is replaced by  full/diagonal matrices $A$ while using  the same reaction. Therefore, the structures of ${\cal L}$ in \reftheo{GenMPMatTKRnew} seems to be necessary for persistence.

\bibliographystyle{plain}

\begin{thebibliography}{10}



\bibitem{Am2} H. Amann, \newblock{Dynamic theory of quasilinear parabolic systems III. Global existence,} {\em  Math Z.} 202 (1989), pp. 219–-250.


\bibitem{BP} Bermon, Abraham and Plemmons, Robert J., {\em  Nonnegative Matrices in the Mathematical Sciences}, Philadelphia: Society for Industrial and Applied Mathematics (1994).

\bibitem{BW} G.J. Butler and P. Waltman. Persistence in dynamical system. {\em Journal of Differential Equations} 63(2):255–263, June 1986.


\bibitem{ET} Eden, A., Foias, C., Nicolaenko, B., and Temam, R. (1994). {\em Exponential Attractors for Evolution Equations}. Research in Applied Maths., Vol. 34, John Wiley and Masson, New York.


\bibitem{JS} O. John and J. Stara, \newblock On the regularity of weak solutions to parabolic systems in two spatial dimensions. \newblock{\em Comm. P.D.E.}, 27(1998), pp. 1159–1170.



\bibitem{GiaS} M. Giaquinta and M. Struwe, \newblock{ On the partial regularity of weak solutions of nonlinear parabolic
systems}. {\em  Math. Z.}, 179(1982),  437--451.








\bibitem{dlejfa} D. Le. \newblock {Existence of Strong and Nontrivial Solutions to Strongly Coupled Elliptic Systems}. {\em J. Funct. Anal.}  272 (2017), no. 11, 4407--4459.



 

\bibitem{dleANS} D. Le, \newblock {Weighted Gagliardo-Nirenberg Inequalities Involving BMO Norms and Solvability of Strongly Coupled Parabolic Systems}. {\em Adv. Nonlinear Stud.} Vol. 16, No. 1(2016), 125--146.

\bibitem{dlebook} D. Le, \newblock{\em Strongly Coupled Parabolic and Elliptic Systems: Existence and Regularity of Strong/Weak Solutions.} De Gruyter, 2018.

\bibitem{dlebook1} D. Le, \newblock{\em Cross Diffusion Systems: Dynamics, Coexistence and Persistence.} De Gruyter, 2022.


\bibitem{dleJMAA} D. Le, \newblock{On the global existence of a generalized Shigesada-Kawasaki-Teramoto system,} {\em J. Math. Anal. App.} 2021.


\bibitem{MAXP} D. Le, \newblock{Some maximum principles for cross difusion systems}, {\em preprint} arXiv:2304.08262.

\bibitem{dleper1} D. Le and T.T. Nguyen, \newblock{Global attractors and uniform persistence for cross diffusion parabolic systems,} {\em Dynamic Systems and Applications} (2007), no. 16, 361--378.

\bibitem{dleper2} D. Le, \newblock{Persistence for a Class of Triangular Cross
	Diffusion Parabolic Systems,} {\em Adv. Nonlinear Stud.} Vol. 5 (2005), 493--514.


\bibitem{HHZ} M. W. Hirsch, H. L. Smith, and X-Q. Zhao. \newblock {Chain transitivity, attractivity, and strong repellors
	for semidynamical systems.} {\em J. Dyn. Diff. Eq.} (1), 13 (1), 107--131 (2001).


\bibitem{lopez} J. López-Gómez and M. Molina-Meyer, \newblock{The maximum principle for cooperative weakly coupled elliptic systems and some applications}. (1994): 383-398.


\bibitem{ST} 
H. L. Smith and H. R. Thieme \newblock{\em Dynamical Systems and Population Persistence.}
American Mathematical Soc., vol. 118, 2011.

\bibitem{Temam} R. Temam, {\em Infinite-Dimensional Dynamical Systems in Mechanics and Physics}, Springer, 1989.


\bibitem{yag} A. Yagi, \newblock{Global solution to some quasilinear parabolic systems in population dynamics}. {\em Nonlin. Anal}. 21
(1993), 603-630.

\bibitem{Z1} E. Zeidler, \newblock{\em Nonlinear Functional Analysis and its Applications, I:Fixed-Point Theorems} Springer-Verlarg New York, 1986.

\end{thebibliography}

\end{document}